\documentclass[journal,twoside]{IEEEtran}

\IEEEoverridecommandlockouts
\usepackage{amsmath,amssymb,amsfonts}
\usepackage{textcomp}
\usepackage{xcolor}
\usepackage{algorithmicx}
\usepackage{algorithm}
\usepackage{algpseudocode}
\usepackage{graphicx}
\usepackage{cite}
\usepackage{color}
\usepackage{url}
\usepackage{array}
\usepackage{multirow}
\usepackage{mathtools}
\usepackage{textcomp}
\usepackage{tikz,amsmath}
 \usetikzlibrary{positioning}
\usetikzlibrary{shapes.geometric,positioning}
\usepackage{siunitx}
\usepackage{circuitikz}
\usepackage{balance}
\usetikzlibrary{shapes.geometric}
\usepackage{soul}
\setstcolor{red}
\setulcolor{red}

\begin{document}

\title{Distributed Control and Optimization of DC Microgrids: A Port-Hamiltonian Approach}

\author{Babak Abdolmaleki,~\IEEEmembership{Member,~IEEE}
and Gilbert Bergna-Diaz,~\IEEEmembership{Member,~IEEE}

\thanks{The authors are with the Department of Electric Power Engineering, Norwegian University of Science and Technology, 7491 Trondheim, Norway (e-mail: babak.abdolmaleki@ntnu.no; gilbert.bergna@ntnu.no).}}

\markboth{
SUBMITTED FOR PUBLICATION. THIS VERSION: August 27, 2021
%IEEE JESTPE
}
{
SUBMITTED FOR PUBLICATION. THIS VERSION: August 27, 2021
%Abdolmaleki \MakeLowercase{\textit{et al.}}: Distributed Control and Optimization of DC Microgrids: A Port-Hamiltonian Approach
}

\maketitle

\begin{abstract}
This article proposes a distributed secondary control scheme that drives a dc microgrid to an equilibrium point where the generators share optimal currents, and their voltages have a weighted average of nominal value. The scheme does not rely on the electric system topology nor its specifications; it guarantees plug-and-play design and functionality of the generators. First, the incremental model of the microgrid system with constant impedance, current, and power devices is shown to admit a port-Hamiltonian (pH) representation, and its passive output is determined. The economic dispatch problem is then solved by the Lagrange multipliers method; the Karush-Kuhn-Tucker conditions and weighted average formation of voltages are then formulated as the control objectives. We propose a control scheme that is based on the Control by Interconnection design philosophy, where the consensus-based controller is viewed as a virtual pH system to be interconnected with the physical one. We prove the regional asymptotic stability of the closed-loop system using Lyapunov and LaSalle theorems. Equilibrium analysis is also conducted based on the concepts of graph theory and economic dispatch. Finally, the effectiveness of the presented scheme for different case studies is validated with a test microgrid system, simulated in both MATLAB/Simulink and OPAL-RT environments.
\end{abstract}
\section{Introduction}
\IEEEPARstart{E}{lectric} power systems are shifting towards the use of more green technologies. To effectively integrate the renewable energy resources, energy storage systems, and electric loads into the power systems, they are interfaced with the grid via power electronic converters and are grouped in the form of microgrids (MGs) easing their control and management\cite{Microgrid}. As a key component of modern power systems, dc microgrids have recently become more attractive\cite{DCReviewQobad}. They are compatible with the dc electric nature of renewable energy resources, energy storage systems, and a majority of electric loads. In addition, compared to ac MGs where control of frequency, phase, reactive power, and power quality are big challenges, control and management of dc grids are inherently simpler\cite{DCReviewQobad}.

In dc MGs, distributed generators (DGs) and loads are connected to the grid via power converters which are either \textit{voltage-controlled} (\textit{grid-forming}) or \textit{current/power-controlled} (\textit{grid-following}).
Grid-forming devices adjust the voltage of their point of common coupling (PCC) to follow a given voltage reference. The grid-following devices, on the other hand, follow some current/power references\cite{Carolina2021}. Therefore, in terms of current/power, the grid-forming DGs and loads are \textit{dispatchable} while the grid-following ones are \textit{non-dispatchable} and can be considered as constant current/power loads. In autonomous MGs, normally a cluster of dispatchable (grid-forming) generators are in charge of shaping the desired voltage level; thus, they should control the dc MG in a collaborative effort.

A common practice to dispatch the current and to  adjust the voltage of grid-forming DGs in a decentralized, communication-free fashion is \textit{droop control}. Despite its simple and robust functionality, this primary controller cannot guarantee desired current-sharing and voltage formation between the DGs\cite{Han2019Review}. To address these shortcomings, the droop characteristic can be corrected by a secondary controller exploiting data exchange either between the DGs and a central control unit or only between the DGs. In the former case, the controller is known as \textit{centralized} secondary control which exhibits a single-point-of-failure to the system and requires a complex communication network between the DGs and the central control unit. Therefore, the distributed control techniques, using neighbor-to-neighbor inter-DG data transmissions, are preferred to the centralized ones\cite{MolzahnReview}.
\subsection{Existing Literature and Research Gap}
Distributed control of dc MGs has already been addressed in many works. A consensus-based proportional current-sharing strategy is proposed in\cite{Nasirian2015} where a dynamic consensus-based estimator is additionally employed to keep the average voltage at nominal value. To reduce the communication burden and to reach faster convergence under this controller, it has respectively  been modified to event-triggered and finite-time versions in\cite{Sahoo2018ET} and \cite{Sahoo2019FT}. In\cite{Peng2020Opt}, a distributed optimal control scheme is proposed under which the DGs can achieve economic current-sharing. Therein, to overcome the initialization and noise robustness problems related to dynamic consensus-based estimation, a modified dynamic consensus-based average voltage observer is used which determines the voltage reference of converters so that the DG currents are shared properly. A somewhat similar control strategy to\cite{Nasirian2015}, but with event-triggered communications, is proposed in\cite{Peng2020ET} under which only information of the DGs' currents are communicated among them. In\cite{Renke2018}, a distributed nonlinear controller is proposed for dc MGs which, instead of droop controller, tunes the DGs voltages so their currents are shared proportionally. To bound the DG voltages within a reasonable range and to guarantee current-sharing among them \textit{to a certain degree}, in\cite{Renke2019}, a containment-consensus-based controller is proposed for dc MGs. A very similar containment-based controller, but with finite-time convergence, is also proposed in\cite{Sahoo2018CFT}.
In the above mentioned works, either the electric network dynamics and electrical system are not taken into account or only a simplified linear algebraic representation of the grid is considered. Consequently, the controller design and system stability may depend on the parameters of the physical system which are subject to modelling uncertainties. 

One way to achieve plug-and-play (PnP) design and operation is to consider the overall system dynamics and to control the system based on energy principles. To do so, in\cite{Cucuzzella2018}, a distributed passivity-based control is proposed for buck-converter-based MGs ensuring proportional current-sharing and average voltage regulation among the DGs. Some similar versions of this controller are presented in\cite{Cucuzzella2019,Trip2018,Trip2019} which demonstrate superior transient system performance. It should be noted that the asymptotic stochastic stability of the controller proposed in\cite{Trip2019} has further been studied in\cite{Silani2020} under varying loads. To reach the desired control objectives in the mentioned works in a finite-time manner, some sliding mode controllers have been developed in\cite{Trip2018SM,Cucuzzella2019Robust}. Moreover, a few consensus-based proportional current-sharing and voltage-balancing controllers, facilitating PnP operations, have been proposed in \cite{Nahata2020ZIP,Nahata2020ZIE} for dc MGs with constant power and exponential loads where the existence and stability of the system equilibria is also studied. Following the same concept and for the sake of PnP functionality of the DGs, in\cite{Sadabadi2021}, a distributed dynamic control strategy is proposed  for voltage balancing and proportional current-sharing among parallel buck converters with the same capacity. The aforementioned works are, however, limited to buck converter-based DGs and \textit{proportional current-sharing} and none of them has considered droop-controlled DGs and their \textit{economic current dispatch}.

\subsection{Contributions}
Motivated by the above literature review, a distributed secondary control strategy for dc MGs with ZIP loads is proposed herein with the following noticeable features. First, in the modeling of the power system, the dynamics of transmission lines and shunt capacitors are considered, loads--and also current/power-controlled DGs-- are modeled by constant-impedance-current-power (ZIP) loads, and the generators are characterized by droop-based grid-forming DGs, encompassing various types of interfacing converters. It is shown that the incremental model of the droop-based MG admits a port-Hamiltonian (pH) representation\cite{Ajan}, and its passive output is defined. Second, drawing inspiration from the Control by Interconnection (CbI) technique  of pH systems\cite{OrtegaCBI}, a distributed consensus-based secondary controller is proposed which drives the MG to an equilibrium point where \textit{i)} the DGs share optimal currents, and \textit{ii)} their weighted-average voltage is the nominal voltage. The voltage weightings are directly related to coefficients of the DGs cost functions and not the electric network and loads. Third, regional asymptotic stability of the system with ZIP loads is demonstrated and it is shown that the system is globally asymptotically stable without the presence of constant power loads (CPLs). Finally, equilibrium analysis is conducted based on the concepts of economic dispatch and graph theory.

The rest of this paper is structured as follows. The MG system modeling and the control aims are formulated in Section~II. Section~III presents the proposed controller and the system stability and equilibrium analyses. The case studies and simulation results are given in Section~IV. Finally, Section~V concludes the paper and discusses future research directions.

Throughout the paper, $\mathbb{R}^{n\times m}$ and $\mathbb{R}^{n}$ stand for the set of $n\times m$ real matrices and $n\times 1$ real vectors, respectively. $\mathrm{diag}\{x_i\}$ indicates a diagonal matrix with $x_i$ being the corresponding diagonal arrays. $\mathrm{col}\{x_i\}$ shows a column vector with the arrays $x_i$. $\mathcal{I}$ is an identity matrix with appropriate dimensions. $\mathbf{0}$ and $\mathbf{1}$ are appropriate all-one and all-zero vectors or matrices. The transpose of a matrix/vector $\mathbf{z}$ is given by $\mathbf{z}^\top$. Given the scalar $x$ or the vector $\mathbf{x}$, $\bar{x}$ and $\bar{\mathbf{x}}$ are their value at the equilibrium point, and $\tilde{x}=x-\bar{x}$ and $\tilde{\mathbf{x}}=\mathbf{x}-\bar{\mathbf{x}}$.

\section{Microgrid Modeling and  Control Objectives}
\subsection{Electric Network, Generators, and ZIP Load Models}
Let $\mathcal{N}_e$, $\mathcal{E}_e$, and $\mathcal{G}_e$, with the cardinalities $n_e^\mathcal{N}$, $n_e^\mathcal{E}$, and $n_e^\mathcal{G}$, be the sets of buses, transmission lines, and grid-forming (voltage-controlled) generators, respectively. Suppose that the transmission lines are modeled by serial resistor-inductor pairs, the buses are modeled by shunt capacitors and ZIP loads, and each generator is modeled by a controllable voltage source which is connected to the grid via a transmission line (See Fig. 1).
\begin{figure}
\centering
\begin{circuitikz}[american,scale=0.69,bigAmp/.style={amp, bipoles/length=1cm}]
\ctikzset{bipoles/length=.69cm}
\scriptsize
        \draw (0,2)
        to[cV, v=$V_i$,] (0,.5) to[short](0,0);
        \draw (0,2)to[short,i=$I_i^{\mathcal{G}_e}$](1,2)       
        to[R=$R_i^{\mathcal{G}_e}$] (2,2) 
        to[L=$L_i^{\mathcal{G}_e}$] (3,2)
        to[short,-*](3.5,2);
        \draw (0,0)to[short](3.5,0) to[C,l=$C_k^{\mathcal{N}_e}$,v<=$V_k^{\mathcal{N}_e}$] (3.5,2)to[short,-*](5,2)to[I,i>_=$I_k^\text{L}$](5,0)to[short](3.5,0);
        \draw (5,2)to[short,i=$I_j^{\mathcal{E}_e}$](6,2)to[R=$R_j^{\mathcal{E}_e}$] (7,2) 
        to[L=$L_j^{\mathcal{E}_e}$] (8,2);
        \draw(6,.5)to[short,i=$ $](5,2);
        \draw (6,.5)to[R] (7,.5) 
        to[L] (8,.5);
        \draw(5,0)to[short](8,0);
        \draw[dotted] (6.25,1.6)to[short,l=$\forall j\in\mathcal{E}_e$](6.25,.9);
        \draw[dotted] (7.75,1.6)to[short](7.75,.9);
        \draw[fill=gray!15] (8,-.2)rectangle(10,2.2)node[midway]{Rest of MG};
        \draw[draw=none,fill=red!5] (-2.2,2.8)rectangle(1.2,4.2);
        \draw[-latex] (-1.4,1.25) -- (-.4,1.25)node[midway,above]{$V_i^\text{ref}$};
        \draw[-latex] (-1,3.5)--(-1.4,3.5)node[midway,below]{$-$};
        \draw[-latex] (-1.5,3.4)--(-1.5,1.35);
        \draw[-latex] (-2.6,1.25) -- (-1.6,1.25)node[near start,above]{$V_\text{nom}$};
        \draw[dashed,blue] (-1.3,.8)rectangle(1.3,1.85);
        \node[blue] at (0,.53) {Converter Dynamics};
        \node[blue] at (0,.2) {\& Internal Controllers};
        \draw (0,3.5) to[bigAmp](-1,3.5);
        \node at (-.35,3.5){$R_i^D$};
        \draw (-1.5,1.25) circle (0.1)node{+};
        \draw (-1.5,3.5) circle (0.1)node{+};
        \draw[-latex] (-1.5,4.2)--(-1.5,3.6)node[near end,left]{$u_i$};
        \draw (0.85,2) ellipse (.07 and .16);
        \draw (0.85,2.16)--(0.85,3.5);
        \draw[-latex] (0.85,3.5)--(0,3.5)node[very near start,above]{$I_i^{\mathcal{G}_e}$};
    \end{circuitikz}
\caption{A droop-based DG connected to a microgrid with ZIP load.}
\end{figure}
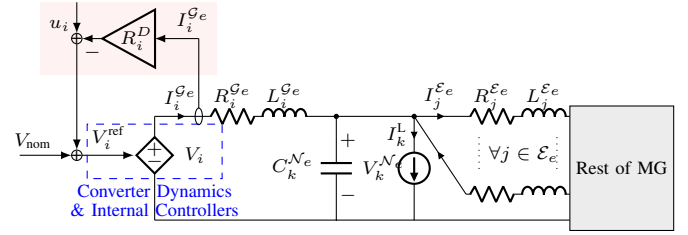

The described electric network can be modeled as two graphs $\mathcal{M}_e$ and $\mathcal{M}_e^\mathcal{G}$ where the buses and transmission lines play the roles of their nodes and edges, respectively. Consider the graph $\mathcal{M}_e=(\mathcal{N}_e,\mathcal{E}_e,\mathcal{B}_e)$ where $\mathcal{N}_e=\{1,\cdots,n_e^\mathcal{N}\}$, $\mathcal{E}_e=\{1,\cdots,n_e^\mathcal{E}\}$, and $\mathcal{B}_e=[b_{kj}]\in \mathbb{R}^{n_e^\mathcal{N}\times n_e^\mathcal{E}}$ are its node set, edge set, and incidence matrix, respectively. Similarly, the graph $\mathcal{M}_e^\mathcal{G}=(\mathcal{N}_e,\mathcal{G}_e,\mathcal{B}_e^\mathcal{G})$ can be defined with the same node set but different edge set $\mathcal{G}_e=\{1,\cdots,n_e^\mathcal{G}\}$ and incidence matrix $\mathcal{B}_e^\mathcal{G}=[b_{ki}^{\mathcal{G}_e}]\in \mathbb{R}^{n_e^\mathcal{N}\times n_e^\mathcal{G}}$. An \textit{incidence matrix} describes the network graph topology by determining the connections between the bus voltages and line currents. For the electric network, one should first consider an \textit{arbitrary} current-flow direction for every line (edge); if current of $j$th line enters node $k$ then $b_{kj}=1$, if it leaves node $k$ then $b_{kj}=-1$, otherwise $b_{kj}=0$. Similarly, if $i$th DG injects current to bus $k$ via an output connector, then $b_{ki}^{\mathcal{G}_e}=1$; otherwise, $b_{ki}^{\mathcal{G}_e}=0$. Note that in this work, the generators are assumed to only inject current to the loads and network and not to absorb it, i.e., $b_{ki}^{\mathcal{G}_e}=-1$ is not considered. 

According to Fig. 1 and based on the system incidence matrices, the dynamics of the droop-based microgrid system are as follows.
\begin{IEEEeqnarray}{rCl}
L_i^{\mathcal{G}_e}\dot{I}_i^{\mathcal{G}_e}&=&V_i-{\sum}_kb_{ki}^{\mathcal{G}_e}V_k^{\mathcal{N}_e}-R_i^{\mathcal{G}_e}I_i^{\mathcal{G}_e},\IEEEyesnumber\IEEEyessubnumber\label{e1a}\\
L_j^{\mathcal{E}_e}\dot{I}_j^{\mathcal{E}_e}&=&-{\sum}_{k} b_{kj}V_k^{\mathcal{N}_e}-R_j^{\mathcal{E}_e}I_j^{\mathcal{E}_e},\IEEEyessubnumber\label{e1b}\\
C_k^{\mathcal{N}_e}\dot{V}_k^{\mathcal{N}_e}&=&{\sum}_{j}b_{kj}I_j^{\mathcal{E}_e}+{\sum}_{i}b^{\mathcal{G}_e}_{ki}I_i^{\mathcal{G}_e}-I_k^\text{L},\IEEEyessubnumber\label{e1c}\\
I_k^\text{L}&=&G_k^\text{cte}V_k^{\mathcal{N}_e}+I_k^\text{cte}+P_k^\text{cte}/V_k^{\mathcal{N}_e},\IEEEyessubnumber\label{e1d}\\
V_i&=&V_i^\text{ref}=V_\text{nom}-R_i^DI_i^{\mathcal{G}_e}+u_i,\IEEEyessubnumber\label{e1e}
\end{IEEEeqnarray}
where $L_i^{\mathcal{G}_e}$, $R_i^{\mathcal{G}_e}$, and $I_i^{\mathcal{G}_e}, \forall i\in\mathcal{G}_e$ are inductance, resistance, and current of $i$th generator transmission line; $L_j^{\mathcal{E}_e}$, $R_j^{\mathcal{E}_e}$, and ${I}_j^{\mathcal{E}_e}, \forall j\in\mathcal{E}_e$ are $j$th line inductance, resistance, and current; $C_k^{\mathcal{N}_e}$, $I_k^\text{L}$, and $V_k^{\mathcal{N}_e},\forall k\in \mathcal{N}_e$ are capacitance, load current, and voltage at bus $k$; $G_k^\text{cte}\geq 0$, $I_k^\text{cte}$, and $P_k^\text{cte}$ are respectively constant conductance, current, and power values of the ZIP load at bus $k$; $V_\text{nom}$, $V_i$, and $V_i^\text{ref}$ are nominal voltage and $i$th DG voltage and its reference value, respectively; $R_i^D$ and $u_i$ are respectively the droop coefficient and correction term (input) of $i$th generator.

There are various types of internal current and/or voltage controllers for converters which are normally designed to be very fast. Hence, in secondary control and optimization design and studies, the following assumptions are usually required.

\textit{Assumption 1:} The grid-forming (voltage-controlled) generators can be modeled as controllable voltage sources so that $V_i=V_i^\text{ref}$. Therefore, considering the well-known droop equation the grid-forming units are characterized by the algebraic relationship $V_i=V_i^\text{ref}=V_\text{nom}-R_i^DI_i^{\mathcal{G}_e}+u_i$.

\textit{Assumption 2:} The grid-following (current-/power- controlled) converters are considered as negative constant current/power loads in the ZIP load model.

Let us define the following global matrices and vectors. $\mathbf{L}_{\mathcal{G}_e}=\mathrm{diag}\{L_i^{\mathcal{G}_e}\}\in\mathbb{R}^{n_e^\mathcal{G}\times n_e^\mathcal{G}}$,
$\mathbf{R}_{\mathcal{G}_e}=\mathrm{diag}\{R_i^{\mathcal{G}_e}\}\in\mathbb{R}^{n_e^\mathcal{G}\times n_e^\mathcal{G}}$, $\mathbf{L}_{\mathcal{E}_e}=\mathrm{diag}\{L_j^{\mathcal{E}_e}\}\in\mathbb{R}^{n_e^\mathcal{E}\times n_e^\mathcal{E}}$, $\mathbf{R}_{\mathcal{E}_e}=\mathrm{diag}\{R_j^{\mathcal{E}_e}\}\in\mathbb{R}^{n_e^\mathcal{E}\times n_e^\mathcal{E}}$, $\mathbf{C}_{\mathcal{N}_e}=\mathrm{diag}\{C_k^{\mathcal{N}_e}\}\in\mathbb{R}^{n_e^\mathcal{N}\times n_e^\mathcal{N}}$, $\mathbf{P}_\text{cte}=\mathrm{col}\{P_k^\text{cte}\}\in\mathbb{R}^{n_e^\mathcal{N}}$, $\mathbf{G}_\text{cte}=\mathrm{diag}\{G_k^\text{cte}\}\in\mathbb{R}^{n_e^\mathcal{N}\times n_e^\mathcal{N}}$, $\mathbf{R}_D=\mathrm{diag}\{R_i^D\}\in\mathbb{R}^{n_e^\mathcal{G}\times n_e^\mathcal{G}}$, $\mathbf{I}_\text{cte}=\mathrm{col}\{I_k^\text{cte}\}\in\mathbb{R}^{n_e^\mathcal{N}}$, $\mathbf{g}_{\mathcal{N}_e}(\mathbf{q}_{\mathcal{N}_e})=\mathrm{diag}\{-C_k^{\mathcal{N}_e}/V_k^{\mathcal{N}_e}\}\in\mathbb{R}^{n_e^\mathcal{N}\times n_e^\mathcal{N}}$, $\mathbf{I}_{\mathcal{G}_e}=\mathrm{col}\{I_i^{\mathcal{G}_e}\}\in\mathbb{R}^{n_e^\mathcal{G}}$, $\mathbf{I}_{\mathcal{E}_e}=\mathrm{col}\{I_j^{\mathcal{E}_e}\}\in\mathbb{R}^{n_e^\mathcal{E}}$, $\mathbf{V}_{\mathcal{N}_e}=\mathrm{col}\{V_k^{\mathcal{N}_e}\}\in\mathbb{R}^{n_e^\mathcal{N}}$, and $\mathbf{u}=\mathrm{col}\{u_i\}\in\mathbb{R}^{n_e^\mathcal{G}}$. Now, with the Hamiltonian $H(\mathbf{x})=0.5\mathbf{x}^\top\mathbf{Q}\mathbf{x}$ where
\begin{IEEEeqnarray}{c}
\mathbf{Q}=\begin{bmatrix}
\mathbf{L}_{\mathcal{G}_e}^{-1} &\mathbf{0}&\mathbf{0}\\
\mathbf{0} & \mathbf{L}_{\mathcal{E}_e}^{-1}&\mathbf{0}\\
\mathbf{0} &\mathbf{0} & \mathbf{C}_{\mathcal{N}_e}^{-1}
\end{bmatrix},\mathbf{x}=\begin{bmatrix}\boldsymbol{\phi}_{\mathcal{G}_e}\\\boldsymbol{\phi}_{\mathcal{E}_e}\\\mathbf{q}_{\mathcal{N}_e}\end{bmatrix}=\mathbf{Q}^{-1}\begin{bmatrix}\mathbf{I}_{\mathcal{G}_e}\\\mathbf{I}_{\mathcal{E}_e}\\\mathbf{V}_{\mathcal{N}_e}\end{bmatrix},\IEEEnonumber
\end{IEEEeqnarray}
one can write the system in the following form.
\begin{IEEEeqnarray}{c}
\Sigma :\begin{cases}
\dot{\mathbf{x}} = \mathbf{F}\nabla H(\mathbf{x})+\mathbf{g}_P(\mathbf{x})\mathbf{P}_\text{cte}+\mathbf{g}\mathbf{u}+\mathbf{E}\\
\mathbf{y}=\mathbf{g}^\top\nabla H(\mathbf{x})
\end{cases},\IEEEyesnumber
\end{IEEEeqnarray}
\begin{IEEEeqnarray}{c}
\mathbf{F} = \mathbf{J}-\mathbf{R}=\begin{bmatrix}
-(\mathbf{R}_{\mathcal{G}_e}+\mathbf{R}_D) &\mathbf{0}&-{\mathcal{B}_e^\mathcal{G}}^\top\\
\mathbf{0} & -\mathbf{R}_{\mathcal{E}_e}&-\mathcal{B}_e^\top\\
\mathcal{B}_e^\mathcal{G}&\mathcal{B}_e & -\mathbf{G}_\text{cte}
\end{bmatrix},\IEEEnonumber\\
\mathbf{g} = \begin{bmatrix} \mathcal{I} \\\mathbf{0}\\\mathbf{0}\end{bmatrix},\mathbf{g}_P(\mathbf{x}) = \begin{bmatrix}\mathbf{0} \\ \mathbf{0} \\ \mathbf{g}_{\mathcal{N}_e}(\mathbf{q}_{\mathcal{N}_e})\end{bmatrix},\mathbf{E} = \begin{bmatrix} \mathbf{1}V_\text{nom} \\\mathbf{0}\\-\mathbf{I}_\text{cte}\end{bmatrix};\IEEEnonumber
\end{IEEEeqnarray}
where $\mathcal{B}_e^\mathcal{G}$ and $\mathcal{B}_e$ are the incidence matrices defined in the preamble of this subsection; $\mathbf{J}=-\mathbf{J}^\top=0.5[\mathbf{F}-\mathbf{F}^\top]$ and $\mathbf{R}=\mathbf{R}^\top=-0.5[\mathbf{F}+\mathbf{F}^\top]$) are the skew-symmetric and symmetric component of $\mathbf{F}$, respectively.

\textit{Assumption 3:} The system $\Sigma$ (2) has a unique equilibrium point $\bar{\mathbf{x}}$. Moreover, $\bar{\mathbf{u}}$ (resp. $\bar{\mathbf{y}}$) are the \textit{equilibrium control} (resp. \textit{equilibrium output}) of (2) at the equilibrium point where
\begin{IEEEeqnarray}{c}
\begin{cases}
\mathbf{0} = \mathbf{F}\nabla H(\bar{\mathbf{x}})+\mathbf{g}_P(\bar{\mathbf{x}})\mathbf{P}_\text{cte}+\mathbf{g}\bar{\mathbf{u}}+\mathbf{E}\\
\bar{\mathbf{y}}=\mathbf{g}^\top\nabla H(\bar{\mathbf{x}})
\end{cases}.\IEEEnonumber
\end{IEEEeqnarray}

The \textit{incremental model} of the system $\Sigma$ for $\tilde{\mathbf{x}}=\mathbf{x}-\bar{\mathbf{x}}$ and $\tilde{\mathbf{u}}=\mathbf{u}-\bar{\mathbf{u}}$ can then be written as the PH system below.
\begin{IEEEeqnarray}{c}
\Tilde{\Sigma} : \begin{cases}
\dot{\tilde{\mathbf{x}}} = [\mathbf{J}-\tilde{\mathbf{R}}(\tilde{\mathbf{x}})]\nabla H(\tilde{\mathbf{x}})+\mathbf{g}\tilde{\mathbf{u}}\\
\tilde{\mathbf{y}}=\mathbf{g}^\top\nabla H(\tilde{\mathbf{x}})
\end{cases},\IEEEyesnumber\label{e8}\\
\tilde{\mathbf{R}}(\tilde{\mathbf{x}}) = \begin{bmatrix}
\mathbf{R}_{\mathcal{G}_e}+\mathbf{R}_D &\mathbf{0}&\mathbf{0}\\
\mathbf{0} & \mathbf{R}_{\mathcal{E}_e}&\mathbf{0}\\
\mathbf{0}&\mathbf{0}& \mathbf{G}_\text{cte}-\mathbf{G}_P(\tilde{\mathbf{q}}_{\mathcal{N}_e})
\end{bmatrix},\IEEEnonumber
\end{IEEEeqnarray}
where $\mathbf{G}_P(\tilde{\mathbf{q}}_{\mathcal{N}_e})=\mathrm{diag}_k\{G_k^P(\tilde{q}_k^{\mathcal{N}_e})\}$ with $G_k^P(\tilde{q}_k^{\mathcal{N}_e})=P_k^\text{cte}(C_k^{\mathcal{N}_e})^2/[\bar{q}_k^{\mathcal{N}_e}(\tilde{q}_k^{\mathcal{N}_e}+\bar{q}_k^{\mathcal{N}_e})]$. 

\textit{Proposition 1:} With the storage function $H(\tilde{\mathbf{x}})=0.5\tilde{\mathbf{x}}^\top\mathbf{Q}\tilde{\mathbf{x}}$, and the passive output $\tilde{\mathbf{y}}$ with respect to $\tilde{\mathbf{u}}$, the system $\Tilde{\Sigma}$ (3) is passive in the following domain.
\begin{IEEEeqnarray}{rCl}
\mathbb{D}&=&\{\tilde{\mathbf{x}}\in\mathbb{R}^{n_e^\mathcal{G}+n_e^\mathcal{E}+n_e^\mathcal{N}}:G_k^\text{cte}>\frac{P_k^\text{cte}(C_k^{\mathcal{N}_e})^2}{\bar{q}_k^{\mathcal{N}_e}(\tilde{q}_k^{\mathcal{N}_e}+\bar{q}_k^{\mathcal{N}_e})}\}.\IEEEyesnumber
\end{IEEEeqnarray}

\textit{Proof:} Since $\mathbf{J}=-\mathbf{J}^\top$, the derivative of the storage function along the trajectories of (3) is
\begin{IEEEeqnarray}{c}
\dot{H}(\tilde{\mathbf{x}})=-(\nabla H(\tilde{\mathbf{x}}))^\top\tilde{\mathbf{R}}(\tilde{\mathbf{x}})\nabla H(\tilde{\mathbf{x}})+\tilde{\mathbf{y}}^\top\tilde{\mathbf{u}} .\IEEEyesnumber
\end{IEEEeqnarray}
On the other hand, the matrix $\tilde{\mathbf{R}}(\tilde{\mathbf{x}})$ is positive definite for all $\tilde{\mathbf{x}}\in\mathbb{D}$. Therefore, the system $\Tilde{\Sigma}$ (3) is passive with the given storage function\cite{Romeo2001}.

\subsection{Economic Dispatch and Near-Nominal Voltage Formation}
Let $\mathcal{C}_i(I_i^{\mathcal{G}_e})=\alpha_i (I_i^{\mathcal{G}_e})^2+\beta_i(I_i^{\mathcal{G}_e})+\gamma_i$ be $i$th generator cost function, where $\alpha_i$, $\beta_i$, and $\gamma_i$ are its parameters. If $I_\text{demand}$ is the total current demand in the power network, then the economic current dispatch problem can be written as the following optimization problem.
\begin{IEEEeqnarray}{c}
\min{\sum}_{i\in\mathcal{G}_e}\mathcal{C}_i(I_i^{\mathcal{G}_e}),\quad\text{s.t.}\quad{\sum}_{i\in\mathcal{G}_e}I_i^{\mathcal{G}_e}=I_\text{demand}.\IEEEnonumber
\end{IEEEeqnarray}
This optimization problem can be solved by Lagrangian method with the following Lagrangian function\cite{Boyd}.
\begin{IEEEeqnarray}{rCl}
L(\mathbf{I}_{\mathcal{G}_e},\lambda)&=&{\sum}_{i\in\mathcal{G}_e}\mathcal{C}_i(I_i^{\mathcal{G}_e})+\lambda(I_\text{demand}-{\sum}_{i\in\mathcal{G}_e}I_i^{\mathcal{G}_e}),\IEEEnonumber
\end{IEEEeqnarray}
where $\lambda$ is \textit{dual variable} or \textit{Lagrange multiplier}. The primal problem is convex; hence, if Slater's condition is satisfied, then the Karush-Kuhn-Tucker (KKT) conditions provide necessary and sufficient conditions for primal-dual optimality of the points as follows\cite{Boyd}.
\begin{IEEEeqnarray}{lCl}
\text{Primal feasibility:}&\quad&\partial L/\partial \lambda=0,\IEEEnonumber\\
\text{Stationary condition:}&\quad&\partial L/\partial I_i^{\mathcal{G}_e}=0,\forall i\in\mathcal{G}_e.\IEEEnonumber
\end{IEEEeqnarray}
This implies that considering a feasible equality constraint in the problem, the KKT optimality conditions are boiled down to the stationary condition\cite{Boyd}
\begin{IEEEeqnarray}{c}
{\lim}_{t\rightarrow \infty}\lambda_i=\lambda_j=\lambda_\text{opt}
\end{IEEEeqnarray}
where $\lambda_i=\partial\mathcal{C}_i/\partial I_i^{\mathcal{G}_e}=2\alpha_iI_i^{\mathcal{G}_e}+\beta_i$ is the incremental cost (Lagrange multiplier) of $i$th DG, and $\lambda_\text{opt}$ is its optimal value. This condition is known as equal incremental costs (EIC) criteria\cite{Wood2013}. Due to the fact that current sharing in power networks depends on the bus-voltage differences and not the absolute values of voltages, theoretically speaking, the above mentioned optimality condition can be satisfied in many voltage levels; i.e., $\lambda_\text{opt}$ can have various values depending on the weighted average of voltages. However, in practice the voltages must be as close as possible to the network's nominal voltage. Therefore, the controller should also guarantee a near-nominal voltage formation which can be formulated as
\begin{IEEEeqnarray}{rCl}
\lim_{t\rightarrow \infty} {\sum}_{i\in\mathcal{G}_e} w_iV_i&=&V_\text{nom}{\sum}_{i\in\mathcal{G}_e} w_i.
\end{IEEEeqnarray}
where $w_i>0,\forall i\in\mathcal{G}_e$ are voltage weightings which are defined later.

\textit{Remark 1:} A special choice of the cost function parameters is $\alpha_i=0.5/I_i^\text{rated}$, $\beta_i=\gamma_i=0$ which turns (6) into the equal current ratios criteria ($I_i^{\mathcal{G}_e}/I_i^\text{rated}=I_j^{\mathcal{G}_e}/I_j^\text{rated}$) underlining the proportional current-sharing, studied in the literature (See e.g., \cite{Renke2018,Renke2019,Sahoo2018CFT,Cucuzzella2018,Cucuzzella2019,Trip2018,Trip2019,Silani2020,Trip2018SM,Cucuzzella2019Robust,Nahata2020ZIP,Nahata2020ZIE,Sadabadi2021}).

\section{Controller Design, Closed-Loop System Equilibrium, and Stability Analysis}
In this section, a distributed controller is proposed for the droop-based MG system to satisfy the control objectives in (6) and (7). The proposed controller relies on both the local and neighborhood measurements of the generators; hence, the generators need to exchange information through a communication network as described next.
\subsection{Communication Network Model}
A communication network between the generators can be modeled as a undirected graph with generators and communication links being its nodes and edges, respectively. Consider the graph $\mathcal{M}_c=(\mathcal{N}_c,\mathcal{E}_c,\mathcal{A})$, where $\mathcal{N}_c=\{1,\cdots,n^\mathcal{N}_c\}$, $\mathcal{E}_c\subseteq \mathcal{N}_c\times\mathcal{N}_c$, and $\mathcal{A}=[a_{ij}]\in \mathbb{R}^{n^\mathcal{N}_c\times n^\mathcal{N}_c}$ are its node set, edge set, and adjacency matrix, respectively. If nodes $i$ and $j$ exchange data, then they are neighbors, $(j,i) \in \mathcal{E}_c$, and $a_{ij}=a_{ji}>0$; otherwise, nodes $i$ and $j$ are not neighbors, $(j,i) \notin \mathcal{E}_c$, and $a_{ij}=a_{ji}=0$. Let $N_i=\{j|(j,i)\in \mathcal{E}_c\}$ and $d_i=\sum_{j\in N_i} a_{ij}$ be neighbor set and in-degree of node $i$, respectively. The Laplacian matrix of $\mathcal{M}_c$ is then $\mathcal{L}=\mathcal{L}^\top\coloneqq\mathcal{D}-\mathcal{A}$, where $\mathcal{D}=\mathrm{diag}\{d_i\}$\cite{Olfati2007}.

\subsection{The Distributed Consensus-Based Control System}
The consensus algorithm\cite{Olfati2007} is an effective technique to perform a distributed solution of the KKT condition in optimization problems (the control objective (6))\cite{MolzahnReview}. Accordingly, we choose the distributed consensus-based integral controller
\begin{IEEEeqnarray}{rCl}
\dot{x}_i^c&=&k_i^I{\sum}_{j\in N_i}a_{ij}(u_j^c-u_i^c),\IEEEyesnumber\IEEEyessubnumber
\end{IEEEeqnarray}
where $u_i^c$ is the data shared between the DGs; $x_i$ is the controller state; $k_i^I>0$ is the integral gain; $a_{ij}$ is the communication weight between DGs $i$ and $j$, defined in the previous subsection. Let us define $\mathbf{x}_c=\mathrm{col}\{x_i^c\}\in \mathbb{R}^{n_e^\mathcal{G}}$, $\mathbf{u}_c=\mathrm{col}\{u_i^c\}\in \mathbb{R}^{n_e^\mathcal{G}}$, and $\mathbf{k}_I=\mathrm{diag}\{k_i\}\in \mathbb{R}^{n_e^\mathcal{G}\times n_e^\mathcal{G}}$. With the Hamiltonian $H_c(\mathbf{x}_c)=0.5\mathbf{x}_c^\top\mathbf{k}_I^{-1}\mathbf{x}_c$, this controller can then be represented as the PH system below.
\begin{IEEEeqnarray}{rCl}
\Sigma_c :\begin{cases}
\dot{\mathbf{x}}_c=\mathbf{g}_c\mathbf{u}_c\\
\mathbf{y}_c=\mathbf{g}_c^\top \nabla H_c(\mathbf{x}_c)
\end{cases},\text{ where } \mathbf{g}_c=-\mathbf{k}_I\mathcal{L}.\IEEEyessubnumber
\end{IEEEeqnarray}
where $\mathcal{L}$ is the Laplacian matrix of the communication network. Now if $\Sigma_c$ (8b) has a feasible equilibrium point, then the incremental model of this linear system can be written as
\begin{IEEEeqnarray}{rCl}
\Tilde{\Sigma}_c :\begin{cases}
\dot{\tilde{\mathbf{x}}}_c=\mathbf{g}_c\tilde{\mathbf{u}}_c\\
\tilde{\mathbf{y}}_c=\mathbf{g}_c^\top \nabla H_c(\tilde{\mathbf{x}}_c)
\end{cases}.\IEEEyessubnumber
\end{IEEEeqnarray}
Therefore, one can write $\dot{H}_c(\tilde{\mathbf{x}}_c)=\tilde{\mathbf{y}}_c^\top\tilde{\mathbf{u}}_c$. Hence, with the storage function $H_c(\tilde{\mathbf{x}}_c)$, the control system $\Tilde{\Sigma}_c$ (8c) is also passive (lossless) with the input $\tilde{\mathbf{u}}_c$ and output $\tilde{\mathbf{y}}_c$.

\subsection{Control by Interconnection of the Incremental Systems}
Now that the incremental model of both physical and control systems are represented as PH systems, one can couple them through the following subsystem.
\begin{IEEEeqnarray}{c}
\Sigma_I :\begin{cases}
\begin{bmatrix}
\mathbf{u}\\\mathbf{u}_c
\end{bmatrix}=\begin{bmatrix}
-\mathbf{r} & -\mathbf{w}^{-1}\\(\mathbf{w}^{-1})^\top & \mathbf{0}
\end{bmatrix}\begin{bmatrix}
\mathbf{y}\\\mathbf{y}_c
\end{bmatrix}+\begin{bmatrix}
\mathbf{b}\\\mathbf{b}_c
\end{bmatrix}
\end{cases},\IEEEyesnumber\IEEEyessubnumber
\end{IEEEeqnarray}
where $\mathbf{b}=\mathrm{col}\{b_i\}$ and $\mathbf{b}_c=\mathrm{col}\{b_i^c\}$ are constant vectors in $\mathbb{R}^{n_e^\mathcal{G}}$; $\mathbf{r}$ and $\mathbf{w}$ are square matrices belonging to $\mathbb{R}^{n_e^\mathcal{G}\times n_e^\mathcal{G}}$. 

\textit{Assumption 4:} The systems $\Sigma$ (2) and $\Sigma_c$ (8b) have feasible equilibrium points which are coupled through the subsystem $\Sigma_I$ (9a) as follows.
\begin{IEEEeqnarray}{c}
\begin{cases}
\begin{bmatrix}
\bar{\mathbf{u}}\\\bar{\mathbf{u}}_c
\end{bmatrix}=\begin{bmatrix}
-\mathbf{r} & -\mathbf{w}^{-1}\\(\mathbf{w}^{-1})^\top & \mathbf{0}
\end{bmatrix}\begin{bmatrix}
\bar{\mathbf{y}}\\\bar{\mathbf{y}}_c
\end{bmatrix}+\begin{bmatrix}
\mathbf{b}\\\mathbf{b}_c
\end{bmatrix}
\end{cases}.\IEEEyessubnumber
\end{IEEEeqnarray}

If \textit{Assumption 4} holds, then the incremental model of $\Sigma_I$ (9a) can be written as the following lossy \textit{interconnection subsystem}\cite{OrtegaCBI}.
\begin{IEEEeqnarray}{c}
\Tilde{\Sigma}_I :\begin{cases}
\begin{bmatrix}
\tilde{\mathbf{u}}\\\tilde{\mathbf{u}}_c
\end{bmatrix}=\begin{bmatrix}
-\mathbf{r} & -\mathbf{w}^{-1}\\(\mathbf{w}^{-1})^\top & \mathbf{0}
\end{bmatrix}\begin{bmatrix}
\tilde{\mathbf{y}}\\\tilde{\mathbf{y}}_c
\end{bmatrix}
\end{cases}.\IEEEyessubnumber
\end{IEEEeqnarray}
Therefore, one has
\begin{IEEEeqnarray}{c}
\tilde{\mathbf{y}}^\top\tilde{\mathbf{u}}+\tilde{\mathbf{y}}_c^\top\tilde{\mathbf{u}}_c=-\tilde{\mathbf{y}}^\top\mathbf{r}\tilde{\mathbf{y}}.\IEEEyessubnumber
\end{IEEEeqnarray}
The configuration used for control by interconnection\cite{OrtegaCBI} of the systems is shown in Fig. 2.
\begin{figure}
\centering
\begin{circuitikz}
\ctikzset{bipoles/length=.7cm}
\normalsize
\draw[fill=gray!50] (0,2)rectangle(1,3.5)node[midway]{$\Sigma$};
\draw[fill=gray!50] (7,2)rectangle(8,3.5)node[midway]{$\Sigma_c$};
\draw[dashed,fill=none] (2.1,2)rectangle(5.5,3.5);\draw[white,fill=white] (4,2) circle (0.26)node{\scriptsize };\node at (4,2){$\Sigma_I$};
\node at (1.2,3){\scriptsize$+$};\draw (1.5,3.2)to[short,i>_=$\:$](1,3.2);\draw(2,3.2)to[R] (3,3.2);\draw[fill=white] (1.7,3.2) circle (0.3)node{\scriptsize };\node at (1.55,3.2){\scriptsize$+$};\node at (1.85,3.2){\scriptsize$-$};\draw(3,3.2)--(3.4,3.2)--(3.4,3.075);\node[diamond,draw,minimum width =0.65cm,minimum height =0.65cm] at (3.4,2.75){};\node at (3.4,2.85){\scriptsize$-$};\node at (3.4,2.65){\scriptsize$+$};\draw(3.4,2.425)--(3.4,2.3)--(1,2.3);\node at (1.2,2.5){\scriptsize$-$};\node at (1.2,2.75){\scriptsize$\mathbf{u}$};\node at (2.5,2.9){\scriptsize$\mathbf{r}$};\node at (1.7,2.75){\scriptsize$\mathbf{b}$};\node at (3.9,2.45){\scriptsize$\mathbf{w}^{-1}\mathbf{y}_c$};\node at (1.2,3.4){\scriptsize$\mathbf{y}$};
%%%%%%%%%%%%%%%%%%%%%%%%%%%%%%%%%%%%%%%%%%%%%%%%%%%%%%%%%%%%%%%%%%%%%%%%%%%%%
\node at (6.8,3){\scriptsize$-$};\draw (7,3.2)to[short,i>_=$\:$](6.5,3.2);\draw (6.5,3.2)--(6.2,3.2)--(6.2,3.075);\draw (6.2,2.775) circle (0.3)node{\scriptsize };\draw[-latex](6.2,3)--(6.2,2.525);\draw (6.2,2.475)--(6.2,2.3)--(7,2.3);\node at (6.8,2.5){\scriptsize$+$};
%%%%%%%%%%%%%%%%%%%%%%%%%%%%%%%%%%
\draw (6.2,3.2)--(5,3.2)--(5,3.075);\node[diamond,draw,minimum width =0.65cm,minimum height =0.65cm] at (5,2.75){};\draw[-latex](5,3)--(5,2.525);\draw (5,2.45)--(5,2.3)--(6.2,2.3);
\node at (6.8,2.75){\scriptsize$\mathbf{y}_c$};\node at (6.8,3.4){\scriptsize$\mathbf{u}_c$};\node at (5.74,2.75){\scriptsize$\mathbf{b}_c$};\node at (4.35,3.1){\scriptsize$(\mathbf{w}^{-1})^\top\mathbf{y}$};
%%%%%%%%%%%%%%%%%%%%%%%%%%%%%%%%%%%%%%%%%%%%%%%%%%%%%%%%%%%%%%%%%%%%%%%%%%%
%%%%%%%%%%%%%%%%%%%%%%%%%%%%%%%%%%%%%%%%%%%%%%%%%%%%%%%%%%%%%%%%%%%%%%%%%%%
\draw[fill=gray!15] (0,0)rectangle(1,1.5)node[midway]{$\Tilde{\Sigma}$};
\draw[fill=gray!15] (7,0)rectangle(8,1.5)node[midway]{$\Tilde{\Sigma}_c$};
\draw[dashed,fill=none] (1.5,0)rectangle(6.6,1.5);\draw[white,fill=white] (4,0) circle (0.26)node{\scriptsize };\node at (4,0){$\Tilde{\Sigma}_I$};
\node at (1.2,1){\scriptsize$+$};\draw (1.5,1.2)to[short,i>_=$\:$](1,1.2);\draw(1.4,1.2)to[R] (2.4,1.2);\draw(2.4,1.2)--(3.4,1.2)--(3.4,1.075);\node[diamond,draw,minimum width =0.65cm,minimum height =0.65cm] at (3.4,0.75){};\node at (3.4,0.85){\scriptsize$-$};\node at (3.4,0.65){\scriptsize$+$};\draw(3.4,0.425)--(3.4,0.3)--(1,0.3);\node at (1.2,0.5){\scriptsize$-$};\node at (1.2,0.75){\scriptsize$\Tilde{\mathbf{u}}$};\node at (1.9,0.9){\scriptsize$\mathbf{r}$};\node at (2.65,0.75){\scriptsize$\mathbf{w}^{-1}\Tilde{\mathbf{y}}_c$};\node at (1.2,1.4){\scriptsize$\Tilde{\mathbf{y}}$};
%%%%%%%%%%%%%%%%%%%%%%%%%%%%%%%%%%%%%%%%%%%%%%%%%%%%%%%%%%%%%%%%%%%%%%%%%%%%%
\node at (6.8,1){\scriptsize$-$};\draw (7,1.2)to[short,i>_=$\:$](6.5,1.2);\node at (6.8,0.5){\scriptsize$+$};
%%%%%%%%%%%%%%%%%%%%%%%%%%%%%%%%%%
\draw (6.5,1.2)--(5,1.2)--(5,1.075);\node[diamond,draw,minimum width =0.65cm,minimum height =0.65cm] at (5,0.75){};\draw[-latex](5,1)--(5,0.525);\draw (5,0.45)--(5,0.3)--(7,0.3);
\node at (6.8,0.75){\scriptsize$\Tilde{\mathbf{y}}_c$};\node at (6.8,1.4){\scriptsize$\Tilde{\mathbf{u}}_c$};\node at (5.9,0.75){\scriptsize$(\mathbf{w}^{-1})^\top\Tilde{\mathbf{y}}$};
\node at (1.2,2){\scriptsize (a)};\node at (1.2,0){\scriptsize (b)};
\end{circuitikz}
\caption{Block (circuit) diagram of the control by interconnection scheme for both non-incremental (a) and incremental (b) system models.}
\end{figure}
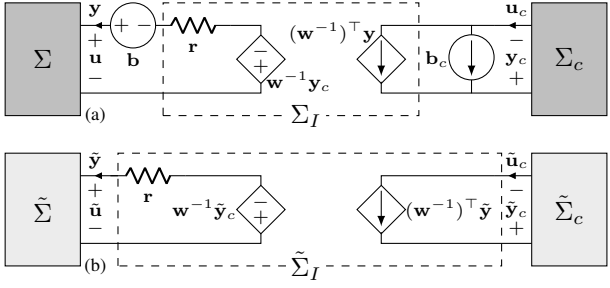

\textit{Proposition 2:} Consider the PH system $\Sigma$ (2) coupled with the controller $\Sigma_c$ (8b) through the interconnection subsystem $\Sigma_I$ (9a) (See Fig. 2). If \textit{Assumption 4} holds and the matrix $\mathbf{R}_D+\mathbf{r}$ is positive-definite, then the equilibrium point of the closed-loop system is \textit{asymptotically stable} in the region
\begin{IEEEeqnarray}{rCl}
\mathbb{S}&=&\{\tilde{\mathbf{x}}_t=[\tilde{\mathbf{x}}^\top,\tilde{\mathbf{x}}_c^\top]^\top:\tilde{\mathbf{x}}\in\mathbb{D},\tilde{\mathbf{x}}_c\in\mathbb{R}^{n_e^\mathcal{G}}\}.\IEEEnonumber
\end{IEEEeqnarray}

\textit{Proof: } Consider the total storage function $H_t(\tilde{\mathbf{x}}_t)=H(\tilde{\mathbf{x}})+H(\tilde{\mathbf{x}}_c)$ for the incremental model of the closed loop-system which has a minimum at the equilibrium point. Taking its derivative and using (3), (8c), and (9d), one has
\begin{IEEEeqnarray}{rCl}
\dot{H}_t(\tilde{\mathbf{x}}_t)&=&\dot{H}(\tilde{\mathbf{x}})+\dot{H}_c(\tilde{\mathbf{x}}_c)=-(\nabla H(\tilde{\mathbf{x}}))^\top\mathbf{T}(\tilde{\mathbf{x}}_t)\nabla H(\tilde{\mathbf{x}});\qquad\IEEEyesnumber\\
\mathbf{T}(\tilde{\mathbf{x}}_t) &=& \begin{bmatrix}
\mathbf{R}_{\mathcal{G}_e}+\mathbf{R}_D+\mathbf{r} &\mathbf{0}&\mathbf{0}\\
\mathbf{0} & \mathbf{R}_{\mathcal{E}_e}&\mathbf{0}\\
\mathbf{0}&\mathbf{0}& \mathbf{G}_\text{cte}-\mathbf{G}_P(\tilde{\mathbf{q}}_{\mathcal{N}_e})
\end{bmatrix}.\IEEEnonumber
\end{IEEEeqnarray}
According to (4), $[\mathbf{G}_\text{cte}-\mathbf{G}_P(\tilde{\mathbf{q}}_{\mathcal{N}_e})]$ is positive-definite for all $\tilde{\mathbf{x}}_t\in\mathbb{S}$, if the closed-loop system has a feasible equilibrium point (\textit{Assumption 4} holds). Moreover, the matrices $\mathbf{R}_{\mathcal{E}_e}$ and $\mathbf{R}_{\mathcal{G}_e}$ are also positive-definite. Therefore, if $\mathbf{R}_D+\mathbf{r}$ is positive-semi-definite, then $\mathbf{T}(\tilde{\mathbf{x}}_t)$ is positive-definite and $\dot{H}_t\leq0, \forall\tilde{\mathbf{x}}_t\in\mathbb{S}$ which proves that the equilibrium point is \textit{stable} in $\mathbb{S}$\cite{Khalil}, with Lyapunov function $H_t$. On the other hand, positive-definiteness of $H_t$ ensures that $\exists \zeta>0$ such that the level set $\Omega_\zeta=\{\tilde{\mathbf{x}}_t\in\mathbb{S}:H_t\leq\zeta\}$ is bounded. Since its requirements are all satisfied, LaSalle's theorem can be applied\cite{Khalil}. According to LaSalle's theorem, every solution starting in $\Omega_\zeta$ converges to the largest invariant set, say $\mathbb{M}$, in $\mathbb{E}=\{\tilde{\mathbf{x}}_t\in\Omega_\zeta:\dot{H}_t=0\}$; i.e., $\tilde{\mathbf{x}}_t\in\mathbb{M}\subseteq\mathbb{E}$ as $t\rightarrow\infty$. Since $\mathbf{T}(\tilde{\mathbf{x}}_t)$ is positive-definite $\forall\tilde{\mathbf{x}}_t\in\mathbb{S}$ and $H(\tilde{\mathbf{x}})$ is quadratic, according to (10), one can write $\mathbb{E}=\{\tilde{\mathbf{x}}_t\in\Omega_\zeta:\tilde{\mathbf{x}}=\mathbf{0},\nabla H(\tilde{\mathbf{x}})=\mathbf{0}\}$ implying that $\dot{\tilde{\mathbf{x}}}=\mathbf{0}$. Therefore, by using (3), (8c), and (9c) it is easy to observe that the motion in this invariant set is governed by $\dot{\tilde{\mathbf{x}}}_c=0, \forall \mathbf{x}_t\in\mathbb{E}$. In other words, the largest invariant set in $\mathbb{E}$ is the equilibrium point; i.e., $\mathbb{M}=\{\tilde{\mathbf{x}}_t\in\Omega_\zeta:\tilde{\mathbf{x}}=\mathbf{0}, \tilde{\mathbf{x}}_c= \mathbf{0}\}$. Therefore, LaSalle's theorem implies \textit{asymptotic stability} of the equilibrium point in $\mathbb{S}$. \hfill $\blacksquare$

\textit{Corollary 1:} Let all the assumptions and conditions of \textit{Propositions 1 and 2} hold. Then, if $\mathbf{P}_\text{cte}= \mathbf{0}$, i.e., if the constant-power loads do not exist in the system, then the equilibrium point of the closed-loop system is \textit{globally asymptotically stable}.

\textit{Proof:} According to (4), since $G_k^\text{cte}\geq 0$ if $P_k^\text{cte}=0$ and all the conditions of \textit{Proposition 1} hold, then one has $\mathbb{D}=\mathbb{R}^{n_e^\mathcal{G}+n_e^\mathcal{E}+n_e^\mathcal{N}}$ and hence $\mathbb{S}=\mathbb{R}^{2n_e^\mathcal{G}+n_e^\mathcal{E}+n_e^\mathcal{N}}$. Moreover, the Lyapunov function $H_t(\tilde{\mathbf{x}}_t)$ is \textit{radially unbounded} as it is in quadratic form. Thus, the equilibrium point is \textit{globally asymptotically stable}, if all the assumptions and conditions of \textit{Proposition 2} hold\cite{Khalil}.\hfill $\blacksquare$

\subsection{Equilibrium (Steady State) Analysis}
\textit{Proposition 3:} Let \textit{Assumption 4} hold. Then, if the communication network is connected, the KKT optimality condition in (6) and the near-nominal voltage formation in (7) with $w_i=\alpha_i^{-1}$ are simultaneously achieved at the equilibrium point of the closed-loop system.

\textit{Proof:} According to (8b), at equilibrium point one has $\mathcal{L}\bar{\mathbf{u}}_c=\mathbf{0}$, where we used the fact that $\mathbf{k}_I$ is positive-definite. If the communication network is connected, then $\mathcal{L}$ has a simple zero eigenvalue\cite{Olfati2007} and therefore $\bar{\mathbf{u}}_c=u_\text{opt}\mathbf{1}$ is the unique solution of $\mathcal{L}\bar{\mathbf{u}}_c=\mathbf{0}$, where $u_\text{opt}$ is the consensus value. Thus, according to (9b) one can write
\begin{IEEEeqnarray}{rCl}
(\mathbf{w}^{-1})^\top\bar{\mathbf{y}}+\mathbf{b}_c&=&u_\text{opt}\mathbf{1}.\IEEEyesnumber\IEEEyessubnumber
\end{IEEEeqnarray}
Let us define $\boldsymbol{\lambda}=\mathrm{col}\{\lambda_i\}\in \mathbb{R}^{n_e^\mathcal{G}}$, $\boldsymbol{\beta}=\mathrm{col}\{\beta_i\}\in \mathbb{R}^{n_e^\mathcal{G}}$ and $\boldsymbol{\alpha}=\mathrm{diag}\{\alpha_i\}\in \mathbb{R}^{n_e^\mathcal{G}\times n_e^\mathcal{G}}$. The KKT condition (6) can then be written as
\begin{IEEEeqnarray}{rCl}
\bar{\boldsymbol{\lambda}}&=&2\boldsymbol{\alpha}\bar{\mathbf{y}}+\boldsymbol{\beta}=\lambda_\text{opt}\mathbf{1}.\IEEEyessubnumber
\end{IEEEeqnarray}
Therefore, if $\mathbf{w}^{-1}=2\boldsymbol{\alpha}$ and $\mathbf{b}_c=\boldsymbol{\beta}$, then one has $u_\text{opt}=\lambda_\text{opt}$ and $\bar{\lambda}_i=\lambda_\text{opt}$. This underlines that the KKT condition is satisfied at the equilibrium point.

Let us further define $\mathbf{V}=\mathrm{col}\{V_i\}\in\mathbb{R}^{n_e^\mathcal{G}}$. From (1e) and (9b) one can write $\bar{\mathbf{V}}=\mathbf{1}V_\text{nom}-\mathbf{R}_D\bar{\mathbf{y}}+\bar{\mathbf{u}}$ and $\bar{\mathbf{u}}=-\mathbf{r}\bar{\mathbf{y}}-\mathbf{w}^{-1}\bar{\mathbf{y}}_c+\mathbf{b}$, and hence $\bar{\mathbf{V}}=\mathbf{1}V_\text{nom}-(\mathbf{R}_D+\mathbf{r})\bar{\mathbf{y}}-\mathbf{w}^{-1}\bar{\mathbf{y}}_c+\mathbf{b}$. Multiplying this equality by $\mathbf{1}^\top\mathbf{w}$ one has
\begin{IEEEeqnarray}{rCl}
\mathbf{1}^\top\mathbf{w}\bar{\mathbf{V}}&=&\mathbf{1}^\top\mathbf{w}\mathbf{1}V_\text{nom}-\mathbf{1}^\top\mathbf{w}[(\mathbf{R}_D+\mathbf{r})\bar{\mathbf{y}}-\mathbf{b}]-\mathbf{1}^\top\bar{\mathbf{y}}_c.\qquad\IEEEnonumber
\end{IEEEeqnarray}
Now if with $k_p\geq 0$ one selects $\mathbf{r}=-\mathbf{R}_D+k_P\mathbf{w}^{-1}\mathcal{L}\mathbf{w}^{-1}$ and $\mathbf{b}=-k_P\mathbf{w}^{-1}\mathcal{L}\boldsymbol{\beta}$, then, using (8b) and the property $\mathbf{1}^\top\mathcal{L}=\mathbf{0}^\top$ of undirected graphs\cite{Olfati2007}, $\mathbf{1}^\top\mathbf{w}\bar{\mathbf{V}}=\mathbf{1}^\top\mathbf{w}\mathbf{1}V_\text{nom}$ can be concluded, which is equivalent to (7) with $w_i=1/(2\alpha_i)$.\hfill $\blacksquare$
\subsection{Implementation of the Proposed Controller}
The proposed controller is presented in matrix form so far. However, in what follows, to better understand its practical implementation, it is formulated in a non-matrix format in terms of the required measurements, parameters, and communication data. If $\mathbf{w}=0.5\boldsymbol{\alpha}^{-1}$, $\mathbf{r}=-\mathbf{R}_D+k_P2\boldsymbol{\alpha}\mathcal{L}2\boldsymbol{\alpha}$, $\mathbf{b}_c=\boldsymbol{\beta}$, and $\mathbf{b}=-k_P2\boldsymbol{\alpha}\mathcal{L}\boldsymbol{\beta}$, then considering (2) and defining $\mathbf{z}_\lambda=\mathrm{col}\{z_i^\lambda\}=-\mathcal{L}\boldsymbol{\lambda}$, and $\mathbf{z}_c=\mathrm{col}\{z_i^c\}=-\mathcal{L}\mathbf{x}_c$, the control system (8b) coupled with the interconnection subsystem (9a) can be written as
\begin{IEEEeqnarray}{c}
\mathbf{u}=\mathbf{R}_D\mathbf{I}_{\mathcal{G}_e}+2\boldsymbol{\alpha}(k_P\mathbf{z}_\lambda-\mathbf{z}_c),\qquad\dot{\mathbf{x}}_c=\mathbf{k}_I\mathbf{z}_\lambda,\IEEEnonumber
\end{IEEEeqnarray}
which can be written in the following scalar format.
\begin{IEEEeqnarray}{c}
\begin{cases}
u_i=R_i^DI_i^{\mathcal{G}_e}+2\alpha_i(k_Pz_i^\lambda-z_i^c)\\
\dot{x}_i^c=k_i^Iz_i^\lambda\\
z_i^\lambda=\sum_{j\in N_i}a_{ij}(\lambda_j-\lambda_i)\\
z_i^c=\sum_{j\in N_i}a_{ij}(x^c_j-x^c_i)\\
\lambda_i=2\alpha_i I_i^{\mathcal{G}_e}+\beta_i
\end{cases}.
\end{IEEEeqnarray}

Fig. 3 depicts a schematic diagram of the proposed controller described in (12). One can see that except for $x_j^c$ and $\lambda_j$, received from the neighboring DGs, the other parameters and variables are locally available for each DG.

\begin{figure}
\centering
\begin{circuitikz}[american,scale=0.69,bigAmp/.style={amp, bipoles/length=1cm}]
\ctikzset{bipoles/length=.69cm}
\scriptsize
        \draw[draw=none,fill=green!5] (-2.2,3.5)rectangle(10.2,9);
        \draw[draw=none,fill=red!5] (-2.2,3.5)rectangle(1.3,4.2);
        \draw[-latex] (-1.5,4.7)--(-1.5,3.5)node[near end,left]{$u_i$};
        \draw (-1.5,4.8) circle (0.1)node{+};
        \draw[-latex] (-1,4.8)--(-1.4,4.8);
        \draw (0,4.8) to[bigAmp](-1,4.8);\node at (-.35,4.8){$R_i^D$};
        \draw[-latex] (.6,7.2)--(0,7.2)node[midway,above]{$z_i^\lambda$};\draw[-latex] (-1,7.2)--(-1.4,7.2);
        \draw (1,6)--(.6,6);
        \draw (0,7.2) to[bigAmp](-1,7.2);\node at (-.35,7.2){$k_P$};
        \draw[-latex] (-1.5,7.1)--(-1.5,6.5);\draw[-latex] (-1.5,5.5)--(-1.5,4.9);
        \draw (-1.5,6.5) to[bigAmp](-1.5,5.5);\node at (-1.5,6.2){$2\alpha_i$};
        \draw (-1.5,7.2) circle (0.1)node{+};
        \draw (1,8.4)--(-1.5,8.4)node[near end,above]{$z_i^c$};
        \draw[-latex] (-1.5,8.4)--(-1.5,7.3)node[very near end,left]{$-$};
        \draw (1,8.1)rectangle(5,8.7)node[midway]{$\sum_{j\in N_i}a_{ij}(x_j^c-x_i^c)$};\draw[-latex](6.15,8.4)--(5,8.4)node[midway,above]{$x_j^c$};
        \draw (1,5.7)rectangle(5,6.3)node[midway]{$\sum_{j\in N_i}a_{ij}(\lambda_j-\lambda_i)$};\draw[-latex](6.15,6)--(5,6)node[midway,above]{$\lambda_j$};
        \draw (1.5,7.2) to[bigAmp](2.5,7.2);\node at (1.75,7.2){$k_i^I$};
        \draw (0.6,6)--(0.6,7.2);\draw[-latex](0.6,7.2)--(1.5,7.2);\draw[-latex](2.5,7.2)--(3.2,7.2);
        \draw (3.2,6.8)rectangle(4,7.6)node[midway]{$\int$};
        \draw[-latex](4,7.2)--(6.15,7.2)node[near start,below]{$x^c_i$};\draw[-latex](4.5,7.2)--(4.5,8.1);
        \draw (0.85,3.5)--(0.85,4.8)node[near start,left]{$I_i^{\mathcal{G}_e}$};
        \draw[-latex] (0.85,4.8)--(0,4.8);
        \draw (1.5,4.8) to[bigAmp](2.5,4.8);\node at (1.85,4.8){$2\alpha_i$};
        \draw[-latex](0.85,4.8)--(1.5,4.8);\draw[-latex](2.5,4.8)--(3.1,4.8);\draw (3.2,4.8) circle (0.1)node{+};\draw[-latex](3.2,4)--(3.2,4.7)node[very near start,left]{$\beta_i$};\draw[-latex](3.3,4.8)--(6.15,4.8)node[near start,below]{$\lambda_i$};\draw[-latex](4,4.8)--(4,5.7);
        \draw[dashed,fill=gray!5] (6.15,3.9)rectangle(9.85,8.9);
        \node at (8,7.4){Neighbor-to-Neighbor};\node at (8,6.4){Inter-DG};\node at (8,5.4){Communication Network};
    \end{circuitikz}
\caption{Schematic diagram of the proposed distributed controller.}
\end{figure}
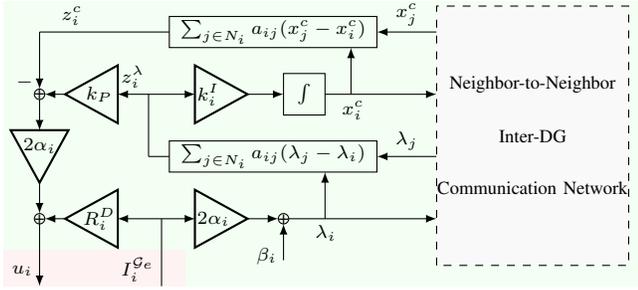

\section{Case Studies and Results}
\begin{figure}
\centering
\begin{circuitikz}[european,scale=0.8,bigAmp/.style={amp, bipoles/length=0.5cm}]
\ctikzset{bipoles/length=.6cm}
\scriptsize
    \draw (0,0)to[R,l=$\mathcal{E}^e_1$,*-*](2,0)to[R,l=$\mathcal{E}^e_2$,-*](4,0)to[R,l=$\mathcal{E}^e_3$,-*](6,0)to[R,l=$\mathcal{E}^e_4$,-*](8,0)to[R,l=$\mathcal{E}^e_5$,-*](10,0);
    \draw (6,-3)to[R,l=$\mathcal{E}^e_8$,*-*](8,-3);
    \draw (6,-3)to[R](4,0);\node at (4.65,-1.7){$\mathcal{E}^e_6$};\draw(8,-3) to[R](10,0);\node at (9.35,-1.7){$\mathcal{E}^e_7$};
    \draw (0,0)to[short](0,-.25)node[sground]{};\draw (2,0)to[short](2,-.25)node[sground]{};\draw (4,0)to[short](4,-.25)node[sground]{};\draw (6,0)to[short](6,-.25)node[sground]{};\draw (8,0)to[short](8,-.25)node[sground]{};\draw (10,0)to[short](10,-.25)node[sground]{};\draw (6,-3)to[short](6,-3.25)node[sground]{};\draw (8,-3)to[short](8,-3.25)node[sground]{};
    \draw (0,1.5)to[R,l=$\mathcal{G}^e_1$](0,0);\draw (2,1.5)to[R,l=$\mathcal{G}^e_2$](2,0);\draw (6,1.5)to[R,l=$\mathcal{G}^e_3$](6,0);\draw (8,1.5)to[R,l=$\mathcal{G}^e_4$](8,0);\draw (6,-1.5)to[R,l=$\mathcal{G}^e_5$](6,-3);\draw (8,-1.5)to[R,a=$\mathcal{G}^e_6$](8,-3);
    \draw[fill=green!6] (0,1.6) circle (0.35)node{DG1};\draw[fill=green!6] (2-.35,1.6-.35) rectangle (2+.35,1.6+.35)node[midway]{DG2};\draw[fill=green!6] (6,1.6) circle (0.35)node{DG3};\draw[fill=green!6] (8-.35,1.6-.35) rectangle (8+.35,1.6+.35)node[midway]{DG4};\draw[fill=green!6] (6,-3+1.6) circle (0.35)node{DG5};\draw[fill=green!6] (8-.35,-3+1.6-.35) rectangle (8+.35,-3+1.6+.35)node[midway]{DG6};
    \node at (0.3,-.2){$\mathcal{N}^e_1$};\node at (2.3,-.2){$\mathcal{N}^e_2$};\node at (6.3,-.2){$\mathcal{N}^e_3$};\node at (8.3,-.2){$\mathcal{N}^e_4$};\node at (6.3,-.2-3){$\mathcal{N}^e_5$};\node at (8.3,-.2-3){$\mathcal{N}^e_6$};\node at (4,.2){$\mathcal{N}^e_7$};\node at (10,.2){$\mathcal{N}^e_8$};
    \draw[dashed,latex-latex,blue] (0.35,1.6) parabola (2-.35,1.6);
    \draw[dashed,latex-latex,blue] (2+.35,1.6)parabola(6-.35,-3+1.6);
    \draw[dashed,latex-latex,blue] (2+.35,1.6)parabola(6-.35,1.6);
    \draw[dashed,latex-latex,blue] (6+.35,1.6)parabola(8-.35,1.6);
    \draw[dashed,latex-latex,blue] (6+.35,-3+1.6)parabola(8-.35,-3+1.6);
    \draw[dashed,latex-latex,blue] (6-.33,1.5) .. controls (5.1,1.3) .. (6-.3,-3+1.8);
    \draw[dashed,latex-latex,blue] (8+.35,1.6) .. controls (10,0) .. (8+.35,-3+1.6);
    \draw[draw=none,fill=red!5] (0,-3.4)rectangle(4.3,-0.9);
    \draw[fill=green!6] (0.5,-1.2) circle (0.15);
    \draw[fill=green!6] (0.35,-1.75) rectangle (0.65,-1.45);
    \draw (0.5,-1.85)node[sground]{};
    \draw (0.2,-2.4) to[R](0.8,-2.4);
    \draw[dashed,latex-latex,blue] (0.1,-2.8) --(0.9,-2.8);
    \draw (0.1,-3.2) --(0.9,-3.2);
    \node [text width=3cm] at (3,-2.17) {\fontsize{7pt}{9.5pt}\selectfont Buck-Based DG\\Boost-Based DG\\Capacitor \& Load\\Transmission Line\\Communication Link\\Electric Connection};
    \draw[draw=none] (1.1,-1.35) rectangle (4,-1.05);
    \draw[draw=none] (1.1,-1.75) rectangle (4,-1.45);
    \draw[draw=none] (1.1,-2.15) rectangle (4,-1.85);
    \draw[draw=none] (1.1,-2.55) rectangle (4,-2.25);
    \draw[draw=none] (1.1,-2.95) rectangle (4,-2.65);
    \draw[draw=none] (1.1,-3.35) rectangle (4,-3.05);
\end{circuitikz}
\caption{Electrical and communication networks of the test microgrid system.\label{TestMG}}
\end{figure}
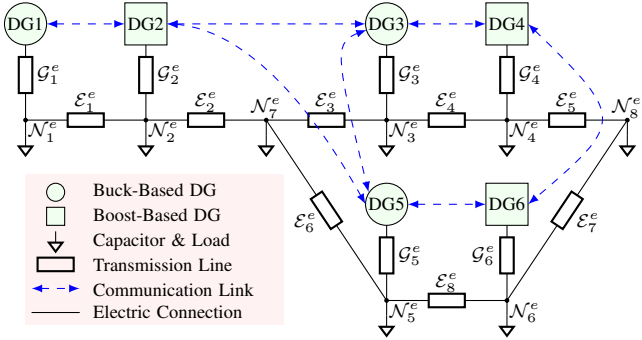

To show the effectiveness of the proposed controller, it is tested on a 48-Volt meshed dc MG, powered by six DGs. The DGs with odd (resp. even) numbers are interfaced to the grid via buck (resp. boost) converters, which are depicted in Fig.~\ref{TestMG} by circles (resp. squares). The electrical and control specifications of the MG shown in Fig.~\ref{TestMG} are given in Table~\ref{TableI}.

\begin{table}
\centering
\caption{The electrical and control specifications of the test MG}
\label{TableI}
\begin{tabular}{|c|c|c|c|c|c|c|c|c|}\hline
\multicolumn{9}{|c|}{DGs' Specifications with Base RL of ($0.5\Omega$,$50\mu H$)}\\ \hline
\multicolumn{3}{|c|}{}&\multicolumn{6}{|c|}{DG Number ($i\in\mathcal{G}_e$)}\\ \cline{4-9}
\multicolumn{3}{|c|}{}	& 1 & 2 & 3 & 4 & 5 & 6 \\ \hline
\multicolumn{3}{|c|}{$I^\text{rated}_i(A)$} & 15 & 6 & 12 & 12 & 10 & 8 \\ \hline
\multicolumn{3}{|c|}{$R_i^D(V/A)$} & 0.2 & 0.5 & 0.25 & 0.25 & 0.3 & 0.375 \\ \hline
\multicolumn{3}{|c|}{$\alpha_i (10^{-1}\$/A^2)$} & 0.8 & 1.9 & 1 & 1.4 & 1.2 & 1.6 \\ \hline
\multicolumn{3}{|c|}{$\beta_i (10^{-1}\$/A)$} & 1 & 2.5 & 1.2 & 1.8 & 1.5 & 2.1 \\ \hline
\multicolumn{3}{|c|}{$\gamma_i(10^{-1}\$)$} & 2 & 5 & 2 & 4 & 3 & 4 \\ \hline
\multicolumn{3}{|c|}{$R_i^{\mathcal{G}_e}$, $L_i^{\mathcal{G}_e}$ (p.u.)} & 0.5 & 0.4 & 0.55 & 0.6 & 0.45 & 0.5\\ \hline
\multicolumn{3}{|c|}{$k_i^P$} &\multicolumn{6}{|c|}{$2$} \\ \hline
\multicolumn{3}{|c|}{$k_i^I$} &\multicolumn{6}{|c|}{$100$} \\ \hline\hline
\multicolumn{9}{|c|}{Line Specifications ($R_i^{\mathcal{E}_e}$, $L_i^{\mathcal{E}_e}$) with Base RL of ($0.5\Omega$,$50\mu H$)}\\ \hline
&\multicolumn{8}{|c|}{Line Number ($j\in\mathcal{E}_e$)}\\ \cline{2-9}
 & 1 & 2 & 3 & 4 & 5 & 6 & 7 &  8 \\ \hline
(p.u.) & 1 &  2 &  2 & 1 & 1 & 3 & 1 & 2 \\ \hline\hline
\multicolumn{9}{|c|}{Bus Specifications}\\ \hline
&\multicolumn{8}{|c|}{Bus Number ($k\in\mathcal{N}_e$)}\\ \cline{2-9}
 & 1 & 2 & 3 & 4 & 5 & 6 & 7 &  8 \\ \hline
$C_k^{\mathcal{N}_e} (F)$ &\multicolumn{8}{|c|}{$22\times10^{-3}$} \\ \hline
$1/G^\text{cte}_k(\Omega)$ & 30 & 20 & 20 & 20 & 30 & 20 & 10 & 10 \\ \hline
$I^\text{cte}_k (A)$ & 0.5 & 0.6 & 0.4 & 0.5 & 0.45 & 0.5 & 0.45 & 0.4 \\ \hline
$P^\text{cte}_k (W)$ &\multicolumn{8}{|c|}{$0.8 G^\text{cte}_k V_n^2$ where $V_n=48 V$} \\ \hline
\end{tabular}
\end{table}

\textit{Remark 2:} According to Assumption 1, to design the secondary controller, the converters are modeled by an equivalent zero-order model as in \eqref{e1e}; thus, the converter dynamics and its internal voltage controller are hidden in Fig.~1 under the dashed blue box. However, in the simulations, Linear Quadratic Regulator (LQR) controller technique is used for the voltage $V_i$ to track its reference $V_i^\text{ref}$\cite{MahdiehDC}. Fig. 5 depicts the converter dynamics and the internal voltage controller. The resistance $R_i$, inductance $L_i$, and capacitance $C_i$ of all the converters are $0.1\Omega$, $2.64mH$, and $2.2mF$, respectively; the input voltage  to the converters $V_i^\text{in}$ of the DGs 1 to 6 are 80, 25, 100, 20, 80, 25 $V$, respectively; $I_i$, $\zeta_i$, and $V_i$ are the states of the system, $m_i$ is the duty cycle given to the PWM generator to produce the switching signal $g_i$ with frequency of 5kHz. To design proper feedback gain matrix $\mathbf{K}_i\in\mathbb{R}^{3\times 3}$, the linearized second-order average model of converters augmented with a voltage-tracker integrator, is used where the output current of the converter capacitor $I_i^{\mathcal{G}_e}$ is considered as an external disturbance, along the lines of\cite{MahdiehDC}.
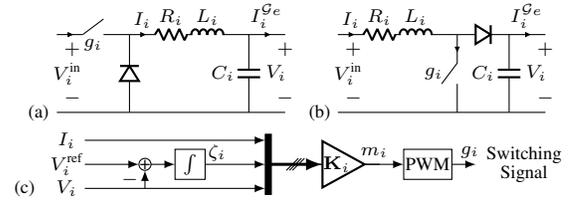
\begin{figure}
\centering
\begin{circuitikz}[american,scale=0.69]
\ctikzset{bipoles/length=.5cm}
\scriptsize
        \draw (.2,1.4)to[open, v=$\:$] (0.2,0.1);\node at(0.1,0.75){$V_i^\text{in}$};
        \draw (-.1,1.5)to[short, i=$\:$](0.4,1.5);\draw(0.4,1.5)--(0.8,1.8);\draw (0.8,1.5)--(1.3,1.5);\draw(1.3,0)to[D=$\:$](1.3,1.5);
        \draw (1.3,1.5)to[short, i=$I_i$](1.8,1.5);\draw (1.8,1.5)to[R=$R_i$](2.4,1.5)--(2.5,1.5);\draw (2.5,1.5)to[L=$L_i$](3.1,1.5)--(3.6,1.5)to[C=$V_i$] (3.6,0);\node at(3.1,0.75){$C_i$};\draw (4.3,1.4)to[open, v=$\:$] (4.3,0.1);\draw(3.6,1.5)to[short,i=$I_i^{\mathcal{G}_e}$](4.3,1.5);\node at(0.6,1.3){$g_i$};
        \draw (-.1,0)--(4.3,0);\node at(-.45,0){(a)};
    \end{circuitikz}
    \begin{circuitikz}[american,scale=0.69]
\ctikzset{bipoles/length=.5cm}
\scriptsize
        \draw (.3,1.4)to[open, v=$\:$] (0.3,0.1);\node at(0.2,0.75){$V_i^\text{in}$};
        \draw (0,1.5)to[short, i=$I_i$](0.5,1.5)to[R=$R_i$](1.1,1.5)--(1.2,1.5)to[L=$L_i$](1.8,1.5)--(2.3,1.5);
        \draw (2.3,1.5)to[short,i=$\:$](2.3,.95)--(2,0.55);\draw (2.3,0.55)--(2.3,0);
        \draw(2.3,1.5)to[D=$\:$](3.3,1.5)to[C=$V_i$] (3.3,0);\node at(2.8,0.75){$C_i$};\draw (4,1.4)to[open, v=$\:$] (4,0.1);\draw(3.3,1.5)to[short,i=$I_i^{\mathcal{G}_e}$](4,1.5);;
        \draw (0,0)--(4,0);\node at(-.35,0){(b)};\node at(1.85,.75){$g_i$};
    \end{circuitikz}\\
    \begin{circuitikz}[european,scale=0.8,bigAmp/.style={amp, bipoles/length=0.8cm}]
\ctikzset{bipoles/length=.5cm}
\scriptsize
    \draw (0.2,0) to[bigAmp](1.4,0);\node at (0.685,0){$\mathbf{K}_i$};\draw[-latex](1.15,0)--(1.8,0)node[near start,above]{$m_i$};\draw (1.8,-.25)rectangle(2.6,0.25);\node at (2.2,0){PWM};\draw[-latex](2.6,0)--(3,0)node[near end,above]{$g_i$};
    \draw[fill=black] (-.5,-.5)rectangle(-.4,.5);\draw[line width=0.4mm,-latex] (-.4,0)--(0.45,0);\draw(-.1,-.1)--(.1,.1);\draw(-.15,-.1)--(.05,.1);\draw(-.05,-.1)--(.15,.1);
    \draw[-latex](-3.5,-.4)--(-.5,-.4);\draw[-latex](-1.5,0)--(-.5,0);\draw[-latex](-3.5,.4)--(-.5,.4);\draw (-2,-.25)rectangle(-1.5,.25);\node at (-1.75,0){$\int$};
    \draw (-2.5,0) circle (0.1)node{+};\draw[-latex](-2.4,0)--(-2,0);\draw[-latex](-2.5,-.4)--(-2.5,-.1)node[midway,left]{$-$};
    \draw[-latex](-3.5,0)--(-2.6,0);\node at (-3.8,0){$V_i^\text{ref}$};\node at (-3.8,-.4){$V_i$};\node at (-3.8,.4){$I_i$};\node at (-1.3,0.2){$\zeta_i$};\node at (-4.5,-.4){(c)};\node at (3.8,0.15){Switching};\node at (3.8,-.15){Signal};
    \end{circuitikz}
\caption{Converter circuit dynamics and internal controller; (a) buck converter, (b) boost converter, and (c) LQR-based voltage controller.}
\end{figure}

\subsection{Controller Performance: Activation and Load Change}
\begin{figure*}
\centering
\includegraphics[width=\textwidth]{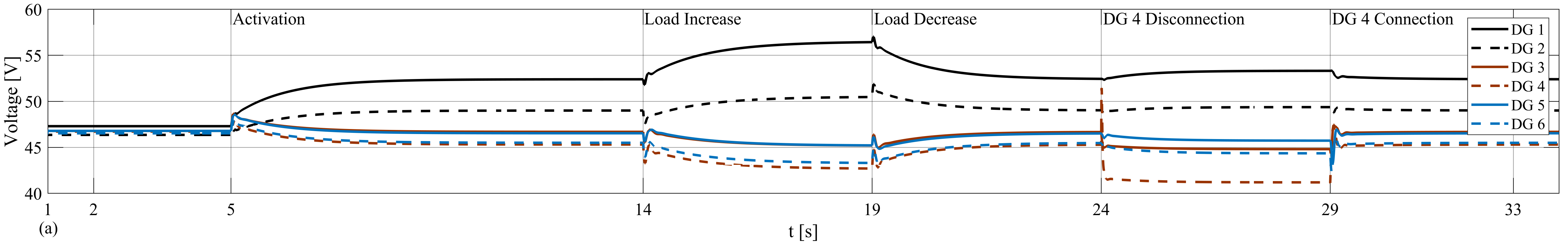}
\includegraphics[width=\textwidth]{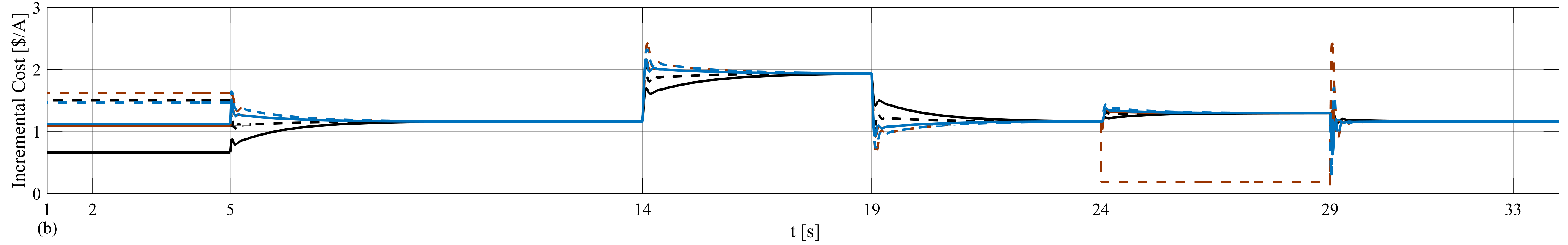}
\includegraphics[width=\textwidth]{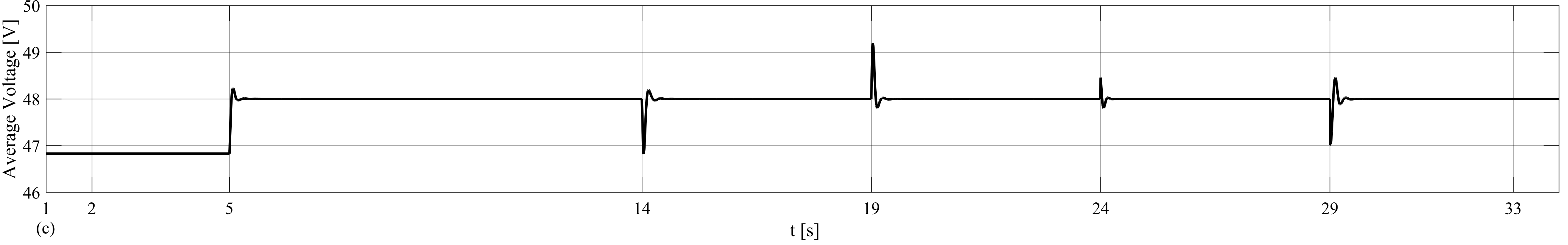}
\caption{Simulation results; (a) the DGs' voltages, (b) the DGs' incremental costs (Lagrange multipliers), and (c) weighted average of the DGs' voltages.\label{Simulation1}}
\end{figure*}
Fig.~\ref{Simulation1} depicts the performance of the MG under the proposed controller in different stages. Before $t=5s$, the MG is operated without the proposed secondary control. Therefore, the DGs voltages are settled away from their nominal voltages so that their average value is deviated from the nominal value $48V$. Moreover, the incremental costs of the DGs have different values which underlines the KKT condition is not satisfied. After activating the controller at $t=5s$, the DGs reach a consensus on their incremental costs and at the same time they form their voltages around the nominal value with a weighted average of nominal voltage. It should be noted that before $t=14s$, only constant impedance and constant current loads are energized. To emphasis the resiliency of the controller, at $t=14s$, the constant power loads at all the buses are activated. One can see that the DGs reach an agreement on a new optimal incremental cost higher than the previous one, which returns to the previous value after deactivating the constant power loads at $t=19s$. It should be emphasized that, over the load change transitions, the average voltage remains unchanged and only transient voltage drifts from the nominal voltage are observed.

\subsection{Controller Performance: Plug-and-Play Ability}
To show the DGs plug-and-play ability under the proposed controller, the 4th DG is disconnected from the grid at $t=24s$ and it is connected back to the grid at $t=29s$. To do so, a corresponding circuit breaker is opened at $t=24s$ to disconnect the DG physically and the communication links related to the DG are all interrupted. Moreover, before closing the breaker at $t=29s$, all the communication links are restored and both sides of the breaker are voltage-synchronized for seamless connection of the DG. According to Fig. \ref{Simulation1}, after disconnecting 4th DG from the grid, other DGs inject more current so they reach consensus on a new optimal incremental cost. Furthermore, one can see that the average voltage of the remaining five DGs still operate at the nominal value while the fourth DG voltage drops to the voltage of the bus number 4. It is also shown that after connecting it back to the grid, the DG immediately participates in the current sharing and voltage formation tasks as before.

\subsection{Real-Time Results From OPAL-RT}
\begin{figure}
\centering
\includegraphics[width=\columnwidth]{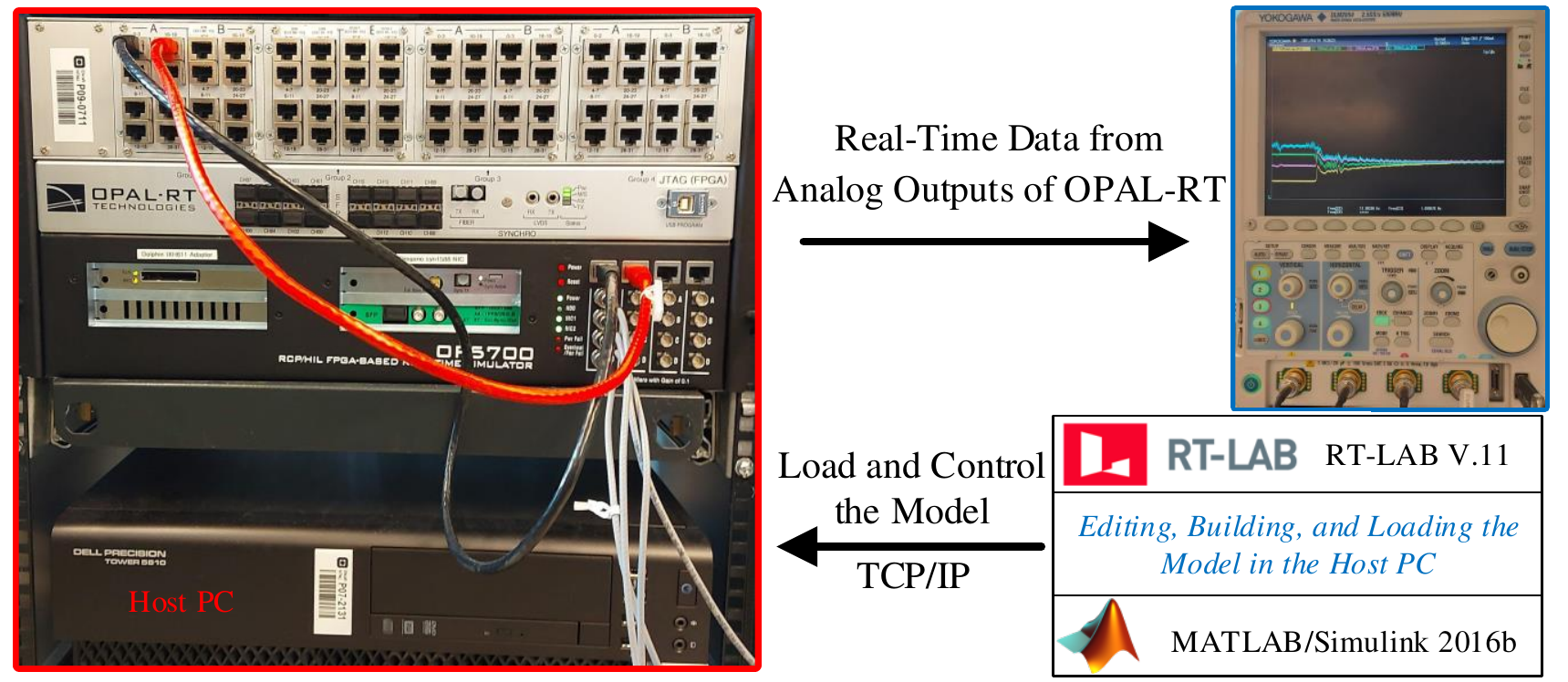}
\caption{Real-time simulation setup.\label{Setup}}
\end{figure}

To verify the real-time effectiveness of the proposed controller, the previous system is built and loaded to an OPAL-RT OP5600 real-time simulator, shown in Fig.~\ref{Setup}. It should be pointed out that, therein, the detailed switching model of the Buck and Boost converters with switching frequency of 5kHz are employed. The selected IGBTs and Diodes have internal resistance of $1m\Omega$ and forward voltage of $0.8V$. The other (passive) components of the converters and their inner voltage controllers are exactly the same as described in the preamble of this Section (Remark 2).

\begin{figure}
\centering
\includegraphics[width=\columnwidth]{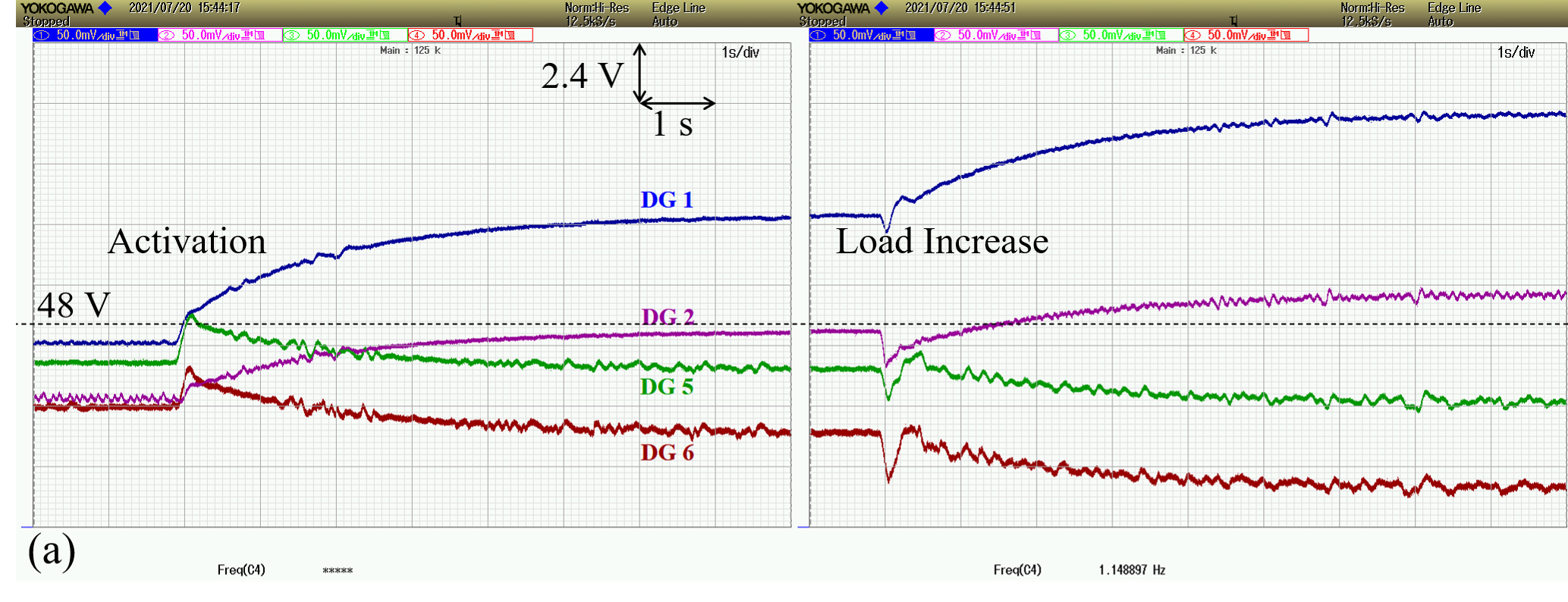}
\includegraphics[width=\columnwidth]{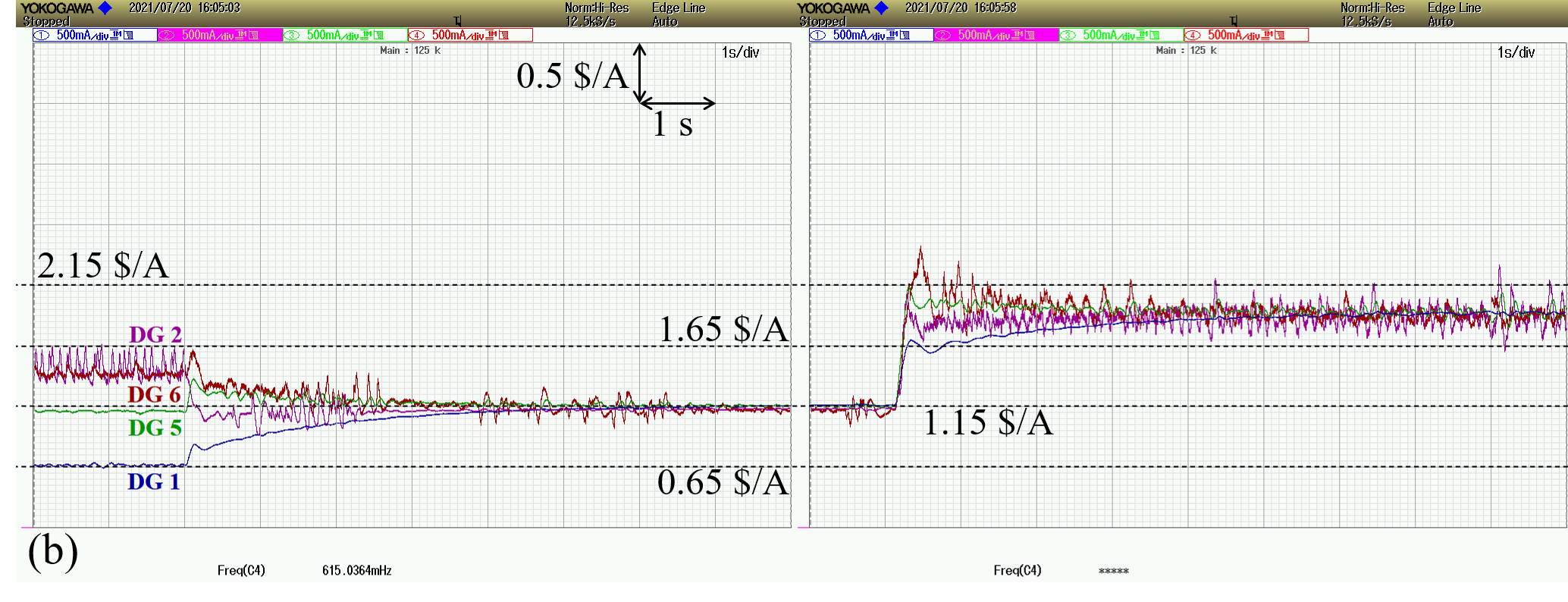}
\includegraphics[width=\columnwidth]{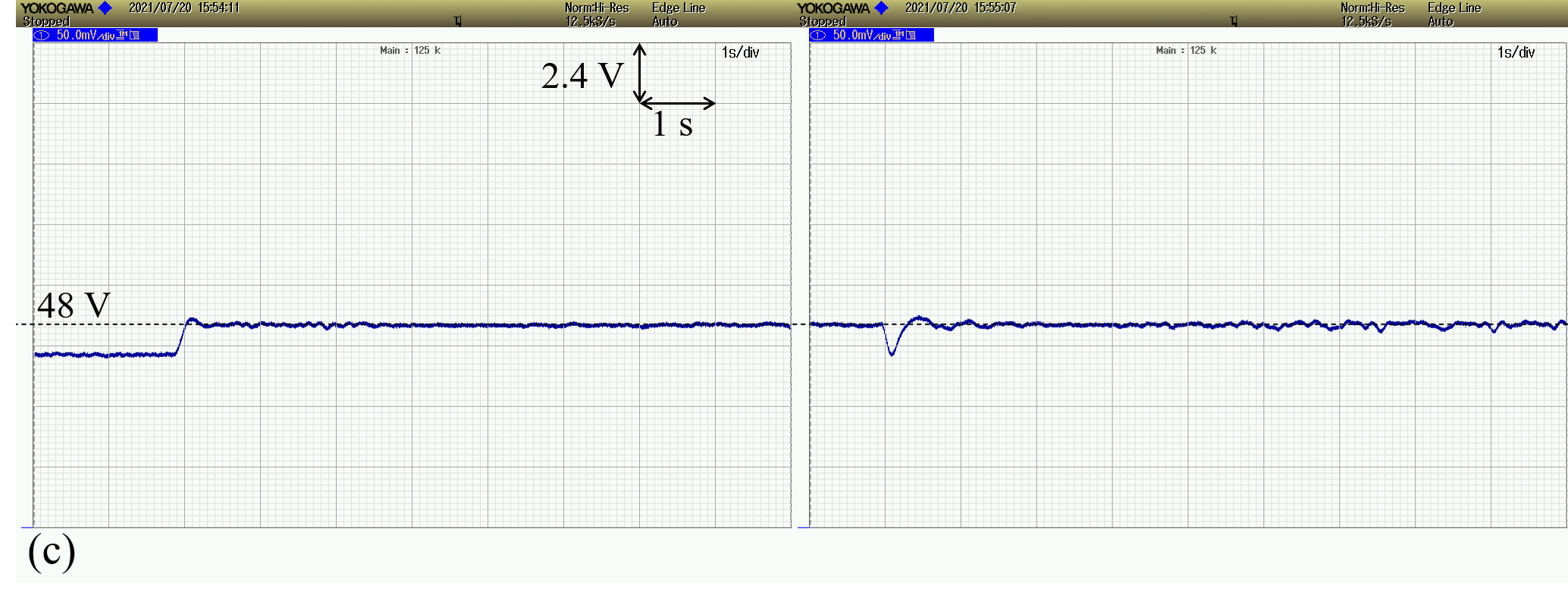}
\caption{Real-time results; (a) the DGs' voltages, (b) the DGs' incremental costs (Lagrange multipliers), and (c) weighted average of the DGs' voltages.\label{RTResults}}
\end{figure}

Fig.~\ref{RTResults} indicates alignment of the real-time system responses with the simulation results in Section IV-A. Due to the input limitation of the oscilloscope only the results for the DGs 1, 2, 5, and 6 are given. After activating the controller, the incremental costs reach a consensus and the voltages reach a formation around the nominal value so that their weighted average settles at the nominal value. The results for the load increase scenario further approves the effectiveness of the proposed control in reaching the current-sharing and voltage-formation control goals, under severe load changes.

\section{Conclusions}
A distributed secondary control technique is proposed for dc MGs with ZIP loads which drives the MG to a point where the KKT optimality condition is satisfied for all the DGs and their weighted average voltage is the nominal value. The closed-loop system (the MG engaged with the proposed controller) is formulated in a port-Hamiltonian representation which is shown to be asymptotically stable by using Lyapunov and LaSalle theorems. It is also shown that the system is globally asymptotically stable without the constant power loads. The effectiveness of the proposed controller for different case studies is verified by adapting it to a test system through both non-real-time and real-time simulations. It should be noted that for the theoretical analyses each DG is modeled by an equivalent zero order model as a controllable voltage source, while, in MATLAB/Simulink simulation and OPAL-RT model the average model and detailed switching model are used, respectively. All in all, the theoretical analyses and case studies demonstrate effectiveness of the proposed controller in achieving the desired control goals.

%\section*{Acknowledgment}

%The authors would like to thank Prabhat Ranjan Bana for his help in obtaining the real-time results from OPAL-RT.

\bibliographystyle{IEEEtran}
\bibliography{IEEEabrv,References}

\begin{thebibliography}{10}
\providecommand{\url}[1]{#1}
\csname url@samestyle\endcsname
\providecommand{\newblock}{\relax}
\providecommand{\bibinfo}[2]{#2}
\providecommand{\BIBentrySTDinterwordspacing}{\spaceskip=0pt\relax}
\providecommand{\BIBentryALTinterwordstretchfactor}{4}
\providecommand{\BIBentryALTinterwordspacing}{\spaceskip=\fontdimen2\font plus
\BIBentryALTinterwordstretchfactor\fontdimen3\font minus
  \fontdimen4\font\relax}
\providecommand{\BIBforeignlanguage}[2]{{%
\expandafter\ifx\csname l@#1\endcsname\relax
\typeout{** WARNING: IEEEtran.bst: No hyphenation pattern has been}%
\typeout{** loaded for the language `#1'. Using the pattern for}%
\typeout{** the default language instead.}%
\else
\language=\csname l@#1\endcsname
\fi
#2}}
\providecommand{\BIBdecl}{\relax}
\BIBdecl

\bibitem{Microgrid}
N.~Hatziargyriou, H.~Asano, R.~Iravani, and C.~Marnay, ``Microgrids,''
  \emph{IEEE Power and Energy Magazine}, vol.~5, no.~4, pp. 78--94, 2007.

\bibitem{DCReviewQobad}
L.~{Meng}, Q.~{Shafiee}, G.~F. {Trecate}, H.~{Karimi}, D.~{Fulwani}, X.~{Lu},
  and J.~M. {Guerrero}, ``Review on control of dc microgrids and multiple
  microgrid clusters,'' \emph{{IEEE} J. Emerg. Sel. Top. Power Electron.},
  vol.~5, no.~3, pp. 928--948, Sep. 2017.

\bibitem{Carolina2021}
C.~Albea-Sanchez, ``{Hybrid dynamical control based on consensus algorithm for
  current sharing in DC-bus microgrids},'' \emph{{Nonlinear Analysis: Hybrid
  Systems}}, vol.~39, Feb. 2021, 100972.

\bibitem{Han2019Review}
Y.~{Han}, X.~{Ning}, P.~{Yang}, and L.~{Xu}, ``Review of power sharing, voltage
  restoration and stabilization techniques in hierarchical controlled dc
  microgrids,'' \emph{IEEE Access}, vol.~7, pp. 149\,202--149\,223, 2019.

\bibitem{MolzahnReview}
D.~K. {Molzahn}, F.~{Dörfler}, H.~{Sandberg}, S.~H. {Low}, S.~{Chakrabarti},
  R.~{Baldick}, and J.~{Lavaei}, ``A survey of distributed optimization and
  control algorithms for electric power systems,'' \emph{{IEEE} Trans. Smart
  Grid}, vol.~8, no.~6, pp. 2941--2962, Nov. 2017.

\bibitem{Nasirian2015}
V.~{Nasirian}, S.~{Moayedi}, A.~{Davoudi}, and F.~L. {Lewis}, ``Distributed
  cooperative control of dc microgrids,'' \emph{{IEEE} Trans. Power Electron.},
  vol.~30, no.~4, pp. 2288--2303, Apr. 2015.

\bibitem{Sahoo2018ET}
S.~Sahoo and S.~Mishra, ``An adaptive event-triggered communication-based
  distributed secondary control for dc microgrids,'' \emph{{IEEE} Trans. Smart
  Grid}, vol.~9, no.~6, pp. 6674--6683, Nov. 2018.

\bibitem{Sahoo2019FT}
------, ``A distributed finite-time secondary average voltage regulation and
  current sharing controller for dc microgrids,'' \emph{{IEEE} Trans. Smart
  Grid}, vol.~10, no.~1, pp. 282--292, Jan. 2019.

\bibitem{Peng2020Opt}
J.~Peng, B.~Fan, and W.~Liu, ``Voltage-based distributed optimal control for
  generation cost minimization and bounded bus voltage regulation in dc
  microgrids,'' \emph{{IEEE} Trans. Smart Grid}, vol.~12, no.~1, pp. 106--116,
  Jan. 2021.

\bibitem{Peng2020ET}
J.~Peng, B.~Fan, Q.~Yang, and W.~Liu, ``Distributed event-triggered control of
  dc microgrids,'' \emph{{IEEE} Syst. J.}, vol.~15, no.~2, pp. 2504--2514, Jun.
  2021.

\bibitem{Renke2018}
R.~Han, L.~Meng, J.~M. Guerrero, and J.~C. Vasquez, ``Distributed nonlinear
  control with event-triggered communication to achieve current-sharing and
  voltage regulation in dc microgrids,'' \emph{{IEEE} Trans. Power Electron.},
  vol.~33, no.~7, pp. 6416--6433, Jul. 2018.

\bibitem{Renke2019}
R.~Han, H.~Wang, Z.~Jin, L.~Meng, and J.~M. Guerrero, ``Compromised controller
  design for current sharing and voltage regulation in dc microgrid,''
  \emph{{IEEE} Trans. Power Electron.}, vol.~34, no.~8, pp. 8045--8061, Aug.
  2019.

\bibitem{Sahoo2018CFT}
S.~Sahoo, D.~Pullaguram, S.~Mishra, J.~Wu, and N.~Senroy, ``A containment based
  distributed finite-time controller for bounded voltage regulation \&
  proportionate current sharing in dc microgrids,'' \emph{Applied Energy}, vol.
  228, pp. 2526--2538, Oct. 2018.

\bibitem{Cucuzzella2018}
M.~Cucuzzella, S.~Trip, and J.~Scherpen, ``A consensus-based controller for dc
  power networks,'' \emph{IFAC-PapersOnLine}, vol.~51, no.~33, pp. 205--210,
  2018, 5th IFAC Conference on Analysis and Control of Chaotic Systems CHAOS
  2018.

\bibitem{Cucuzzella2019}
M.~Cucuzzella, K.~C. Kosaraju, and J.~M.~A. Scherpen, ``Distributed
  passivity-based control of dc microgrids,'' in \emph{{Proc.} American Control
  Conference (ACC)}, Philadelphia, PA, USA, Jul. 2019, pp. 652--657.

\bibitem{Trip2018}
S.~Trip, R.~Han, M.~Cucuzzella, X.~Cheng, J.~Scherpen, and J.~Guerrero,
  ``Distributed averaging control for voltage regulation and current sharing in
  dc microgrids: Modelling and experimental validation,''
  \emph{IFAC-PapersOnLine}, vol.~51, no.~23, pp. 242--247, 2018, 7th IFAC
  Workshop on Distributed Estimation and Control in Networked Systems NECSYS
  2018.

\bibitem{Trip2019}
S.~{Trip}, M.~{Cucuzzella}, X.~{Cheng}, and J.~{Scherpen}, ``Distributed
  averaging control for voltage regulation and current sharing in dc
  microgrids,'' \emph{{IEEE} Control Syst. Lett.}, vol.~3, no.~1, pp. 174--179,
  Jan. 2019.

\bibitem{Silani2020}
A.~Silani, M.~Cucuzzella, J.~M.~A. Scherpen, and M.~J. Yazdanpanah, ``Passivity
  properties for regulation of dc networks with stochastic load demand,'' in
  \emph{{Proc.} 21st IFAC World Congress}, Berlin, Germany, Jul. 2020.

\bibitem{Trip2018SM}
S.~Trip, M.~Cucuzzella, C.~D. Persis, X.~Cheng, and A.~Ferrara, ``Sliding modes
  for voltage regulation and current sharing in dc microgrids,'' in
  \emph{{Proc.} American Control Conference (ACC)}, Milwaukee, WI, USA, Jun.
  2018, pp. 6778--6783.

\bibitem{Cucuzzella2019Robust}
M.~{Cucuzzella}, S.~{Trip}, C.~{De Persis}, X.~{Cheng}, A.~{Ferrara}, and
  A.~{van der Schaft}, ``A robust consensus algorithm for current sharing and
  voltage regulation in dc microgrids,'' \emph{{IEEE} Trans. Control Syst.
  Technol.}, vol.~27, no.~4, pp. 1583--1595, Jul. 2019.

\bibitem{Nahata2020ZIP}
P.~{Nahata} and G.~{Ferrari-Trecate}, ``On existence of equilibria, voltage
  balancing, and current sharing in consensus-based dc microgrids,'' in
  \emph{{Proc.} European Control Conference (ECC)}, St. Petersburg, Russia, May
  2020, pp. 1216--1223.

\bibitem{Nahata2020ZIE}
P.~Nahata, M.~S. Turan, and G.~Ferrari-Trecate, ``Consensus-based current
  sharing and voltage balancing in dc microgrids with exponential loads,''
  \emph{arXiv preprint arXiv:2007.10134}, 2020.

\bibitem{Sadabadi2021}
M.~S. Sadabadi, ``A distributed control strategy for parallel dc-dc
  converters,'' \emph{{IEEE} Control Syst. Lett.}, vol.~5, no.~4, pp.
  1231--1236, Oct. 2021.

\bibitem{Ajan}
A.~Van Der~Schaft and D.~Jeltsema, \emph{Port-Hamiltonian Systems Theory: An
  Introductory Overview}.\hskip 1em plus 0.5em minus 0.4em\relax Now
  Foundations and Trends, 2014.

\bibitem{OrtegaCBI}
R.~Ortega, A.~van~der Schaft, F.~Castanos, and A.~Astolfi, ``Control by
  interconnection and standard passivity-based control of port-hamiltonian
  systems,'' \emph{{IEEE} Trans. Autom. Control}, vol.~53, no.~11, pp.
  2527--2542, Dec. 2008.

\bibitem{Romeo2001}
R.~Ortega, A.~Van Der~Schaft, I.~Mareels, and B.~Maschke, ``Putting energy back
  in control,'' \emph{IEEE Control Systems Magazine}, vol.~21, no.~2, pp.
  18--33, 2001.

\bibitem{Boyd}
S.~Boyd and L.~Vandenberghe, \emph{Convex Optimization}.\hskip 1em plus 0.5em
  minus 0.4em\relax USA: Cambridge University Press, 2004.

\bibitem{Wood2013}
A.~J. Wood, B.~F. Wollenberg, and G.~B. Shebl{\'e}, \emph{Power Generation,
  Operation, and Control}, 3rd~ed.\hskip 1em plus 0.5em minus 0.4em\relax
  Wiley-Interscience, 2013.

\bibitem{Olfati2007}
R.~Olfati-Saber, J.~A. Fax, and R.~M. Murray, ``Consensus and cooperation in
  networked multi-agent systems,'' \emph{Proceedings of the IEEE}, vol.~95,
  no.~1, pp. 215--233, Jan. 2007.

\bibitem{Khalil}
H.~K. Khalil, \emph{Nonlinear Systems}, 3rd~ed.\hskip 1em plus 0.5em minus
  0.4em\relax Englewood Cliffs, NJ, USA: Prentice Hall, 2002.

\bibitem{MahdiehDC}
M.~S. {Sadabadi}, Q.~{Shafiee}, and A.~{Karimi}, ``Plug-and-play robust voltage
  control of dc microgrids,'' \emph{{IEEE} Trans. Smart Grid}, vol.~9, no.~6,
  pp. 6886--6896, Nov. 2018.

\end{thebibliography}

\end{document}